\documentclass[11pt]{article}
\usepackage{mathrsfs}
\usepackage{algorithm}
\usepackage{algpseudocode}
\usepackage{amsmath}
\usepackage{amssymb}
\usepackage{bm}
\usepackage{booktabs}
\usepackage{epsfig}
\usepackage{epstopdf}
\usepackage{flafter}
\usepackage{float}
\usepackage{ifthen}
\usepackage{multirow}
\usepackage{amsfonts}
\usepackage{makecell} 
\usepackage{pifont}
\usepackage{subfigure}
\usepackage{times}
\usepackage{threeparttable}
\usepackage{url} 
\usepackage{subfigure,graphicx}
\usepackage{float}
\usepackage{color,multirow}
\usepackage{makecell}
\usepackage{setspace}
\usepackage[nocompress]{cite} 
\usepackage{url} 
\usepackage{ifthen} 
\usepackage[T1]{fontenc}
\usepackage{graphicx}
\usepackage{changepage}
\usepackage[colorlinks,linkcolor=red,anchorcolor=gray,citecolor=red,urlcolor=blue]{hyperref}
\usepackage{graphics}
\usepackage{booktabs}

\usepackage{indentfirst}

\def\sqr#1#2{{\vcenter{\vbox{\hrule height.#2pt
\hbox{\vrule width.#2pt height#1pt \kern#1pt \vrule width.#2pt}
\hrule height.#2pt}}}}
\def\signed #1{{\unskip\nobreak\hfil\penalty50
\hskip2em\hbox{}\nobreak\hfil#1
\parfillskip=0pt \finalhyphendemerits=0 \par}}
\def\endpf{\signed {$\sqr69$}}

\def\dbR{{\mathop{\rm l\negthinspace R}}}
\def\dbC{{\mathop{\rm l\negthinspace\negthinspace\negthinspace C}}}

\def\3n{\negthinspace \negthinspace \negthinspace }
\def\2n{\negthinspace \negthinspace }
\def\1n{\negthinspace }

 \def\sA{\mathscr{A}}  \def\cA{{\cal A}}  
     
\def\dbC{\mathbb{C}}     
     
\def\dbE{\mathbb{E}}   \def\cE{{\cal E}}  
\def\dbF{\mathbb{F}}   \def\cF{{\cal F}}  
     
   \def\cH{{\cal H}}  
     
   \def\cJ{{\cal J}}  
     
   \def\cL{{\cal L}}  
     
\def\dbN{\mathbb{N}}     
   \def\cO{{\cal O}}  
\def\dbP{\mathbb{P}}     
     
\def\dbR{\mathbb{R}}     
\def\dbS{\mathbb{S}}

    \def\BY{{\bf Y}} \def\By{{\bf y}}

\def\ds{\displaystyle}

\def\={\buildrel \triangle \over =}

\def\lg{\langle}
\def\rg{\rangle}
\def\a{\alpha}

\def\g{\gamma}

\def\vep{\varepsilon}

\def\k{\kappa}
\def\l{\lambda}

\def\si{\sigma}
\def\t{\times}
\def\th{\theta}

\def\i{\infty}
\def\vp{\varphi}
\def\ns{\noalign{\ss} }

\def\G{\Gamma}
\def\D{\Delta}

\def\Si{\Sigma}

\def\O{\Omega}

\def\no{\noindent}

\def\ss{\smallskip}
\def\ms{\medskip}
\def\bs{\bigskip}
\def\q{\quad}
\def\qq{\qquad}

\def\ol{\overline}
\def\tl{\tilde}
\def\ch{\check}
\def\max{\mathop{\rm max}}
\def\min{\mathop{\rm min}}

\def\sup{\mathop{\rm sup}}

\def\h{\widehat}
\def\wt{\widetilde}

\def\dist{\hbox{\rm dist$\,$}}

\def\supp{\hbox{\rm supp$\,$}}

\def\({\Big (}
\def\){\Big )}
\def\[{\Big[}
\def\]{\Big]}
\def\bel{\begin{equation}\label}
\def\ee{\end{equation}}
\def\bt{\begin{theorem}}
\def\bcd{\begin{condition}}
\def\ecd{\end{condition}}
\def\et{\end{theorem}}
\def\bc{\begin{corollary}}
\def\ec{\end{corollary}}
\def\bde{\begin{definition}}
\def\ede{\end{definition}}
\def\bl{\begin{lemma}}
\def\el{\end{lemma}}
\def\bp{\begin{proposition}}
\def\ep{\end{proposition}}
\def\br{\begin{remark}}
\def\er{\end{remark}}
\def\ba{\begin{array}}
\def\ea{\end{array}}
\def\ns{\noalign{\ms}}
\def\ds{\displaystyle}

\def\square#1{\vbox{\hrule\hbox{\vrule height#1%
\kern#1\vrule}\hrule}}
\def\rectangle#1#2{\vbox{\hrule\hbox{\vrule height#1%
\kern#2\vrule}\hrule}}

\font\tenbb=msbm10 \font\sevenbb=msbm7 \font\fivebb=msbm5

\newfam\bbfam
\scriptscriptfont\bbfam=\fivebb \textfont\bbfam=\tenbb
\scriptfont\bbfam=\sevenbb

\newtheorem{lemma}{Lemma}[section]
\newtheorem{remark}{Remark}[section]

\newtheorem{theorem}{Theorem}[section]
\newtheorem{corollary}{Corollary}[section]

\newtheorem{definition}{Definition}[section]
\newtheorem{proposition}{Proposition}[section]
\newtheorem{condition}{Condition}[section]

\makeatletter

\@addtoreset{equation}{section}
\makeatother

\allowdisplaybreaks
\everymath{\displaystyle}

\begin{document}
\title{\bf Exact Controllability for a Refined stochastic hyperbolic Equation with  Internal Controls\thanks{This work is  supported by the NSF of China under grants 12025105 and  12401586. }}
\author{Zengyu Li \thanks{School of Mathematics, Sichuan
University, Chengdu, 610064, China.  {\small\it
E-mail:} {\small\tt lizengyu@stu.scu.edu.cn}. }, \q
Zhonghua Liao\thanks{School of Mathematics, Sichuan
University, Chengdu, 610064, China.  {\small\it
E-mail:} {\small\tt zhonghualiao@yeah.net}. } \q
and \q  Qi L\"{u}\thanks{School of Mathematics, Sichuan University, Chengdu, 610064, China. {\small\it E-mail:} {\small\tt lu@scu.edu.cn}. }}
\date{}
\maketitle
\begin{abstract}
We establish the internal exact controllability of a refined stochastic hyperbolic equation by deriving a suitable observability inequality via Carleman estimates for the associated backward stochastic hyperbolic equation. In contrast to existing results on boundary exact controllability--which require  longer waiting times, we demonstrate that the required waiting time for internal exact controllability in stochastic hyperbolic equations coincides exactly with that of their deterministic counterparts. 
\end{abstract}

\no{\bf 2020 Mathematics Subject Classification}. 93B05, 93B07.

\bs

\no{\bf Key Words}. Stochastic hyperbolic equation, internal exact controllability, observability estimate, Carleman estimate.

\section{Introduction}\label{sec1}

Let $T>0$, and ($\Omega,\cF,\mathbf{F},\dbP$) be a complete filtered probability space on which a one dimensional Brownian motion $W(\cdot)$ is defined and $\mathbf{F}=\{\cF_{t}\}_{t\in[0,T]}$ is the natural filtration generated by $W(\cdot)$. Denote by $\dbF$ the progressive $\sigma$-field with respect to $\mathbf{F}$.

Let $H$ be a Banach space. Denote by
$L^{2}_{\dbF}(0, T; H)$ the Banach space consisting of all
$H$-valued $\dbF$-adapted process $X(\cdot)$ such that
$\dbE(\|X(\cdot)\|_{L^2(0, T; H)}^2)< +\infty$,
by $L^{\infty}_{\dbF}(0, T; H)$ the Banach space of all $H$-valued
$\dbF$-adapted essentially bounded processes,
by $L^\infty_{\dbF}(\O;W^{1,\i}(0,T;H))$ the Banach space of all $H$-valued $\dbF$-adapted essentially bounded processes with their first order weak derivatives being still essentially bounded, by $L^2_{\dbF}(\O;H_0^{1}(0,T;H))$ the Banach space of all $H$-valued $\dbF$-adapted  processes $u$ such that their first order weak derivatives belong to $L^{2}_{\dbF}(0, T; H)$ and $u(0)=u(T)=0$, $\dbP$-a.s., 
and by $L^{2}_{\dbF}(\Omega;C([0,T];H))$ the Banach space of all $H$-valued $\dbF$-adapted continuous processes $X(\cdot)$ such that $\dbE(\|X(\cdot)\|_{C([0,T];H)}^{2})<+\infty$. All of the above spaces are equipped with their canonical norms.

Let $G\subset\dbR^{n}$ (where $n\in\dbN$) be a bounded domain with a $C^{2}$ boundary $\Gamma$. Let
$
Q\=(0,T)\times G
$ and $\Si\=(0,T)\times \G$. Let $(a^{jk}(\cdot))_{1\leq j,k\leq n}\in C^{2}(\overline{G};\dbR^{n\t n})$ satisfy  that\vspace{-2mm}
\begin{equation}\label{1.1}
a^{jk}(x)=a^{kj}(x),\q \forall x\in\ol{G},\q j,k=1,\cdots,n,
\end{equation}
and\vspace{-2mm}
\begin{equation}\label{1.2}
\sum_{j,k=1}^{n}a^{jk}(x)\xi_{j}\xi_{k}\geq \a_{0}|\xi|^{2},\q
\forall(x,\xi)=(x,\xi_{1},\cdots,\xi_{n})\in\ol{G}\t\dbR^{n},
\end{equation}
for some $\a_0>0$.

Consider the following controlled refined stochastic hyperbolic equation:\vspace{-2mm}
\begin{equation}\label{1.3}
\begin{cases}
dy=(\hat{y}+a_{5}f)dt+(a_{3}y+f)dW(t) & \text{in}\q Q,\\
d\hat{y}-\sum_{j,k=1}^{n}(a^{jk}y_{x_{j}})_{x_{k}}dt=(a_1y+a_{4}g+\chi_{G_0}h)dt+(a_{2}y+g)dW(t) &\text{in}\quad Q,\\
y=0 & \text{on} \quad \Sigma,\\
y(0)=y_{0},\quad \hat{y}(0)=\hat{y}_{0} & \text{in}\quad G,
\end{cases}
\end{equation}
where $(y_{0},\hat{y}_{0})\in H^{1}_{0}(G)\t L^{2}(G)$, with $a_{1},a_{2},a_{4}\in L^{\infty}_{\dbF}(0,T;L^{\infty}(G))$,  $a_{3},a_{5}\in L_{\dbF}^{\infty}(\Omega;W^{1,\infty}(0,T;$ $L^{\infty}(G)))$, and $f\in L_{\dbF}^{2}(0,T;H^{1}_{0}(G))$, $g\in L_{\dbF}^{2}(0,T;L^{2}(G))$, $h\in L_{\dbF}^{2}(0,T;L^{2}(G))$ are three controls and $G_0\subset G$ is a suitable subset defined in \eqref{1.11}. By the classical wellposedness result for stochastic evolution equations, we know that equation \eqref{1.3} admits a unique weak solution $(y,\hat{y})\in L^2_{\dbF}(\O;C([0,T];H^{1}_{0}(G)))\t L^2_{\dbF}(\O;C([0,T];L^{2}(G)))$.

\begin{remark}
A distinctive feature of equation \eqref{1.3} lies in the inclusion of the terms $a_{5}fdt$ and $a_{4}gdt$. As observed in \cite{LZ21},  incorporating controls into the diffusion terms (i.e., $f dW(t)$ and $gdW(t)$) inevitably induces secondary effects on the drift terms in realistic models. This inherent property of stochastic systems substantially complicates the analysis. For simplicity, in this work, we assume that these secondary effects depend linearly on $f$ and $g$.  
\end{remark}

The classical stochastic hyperbolic equation is introduced to model the vibration of strings and membranes perturbed by random forces, as well as the propagation of waves in a random environment (see, e.g., \cite[Chapter 2]{DKMNX}, \cite{Funaki83}). However, as shown in \cite[Chapter 10]{LZ21}, this equation fails to be exactly controllable even when controls are applied to both the drift and diffusion terms throughout the entire domain.  This differs significantly from the controllability properties of deterministic wave equations. Exact controllability property  is a fundamental property of control systems. It plays a key role in analyzing system stabilizability, achieving disturbance decoupling, establishing the existence of optimal controls, and determining the nontriviality of necessary conditions in optimal control problems with endpoint constraints. As a result, when establishing mathematical models for control systems, in addition to the accuracy of the model, whether the system is exactly controllable also serves as a crucial criterion. 
Inspired by this and guided by principles from statistical mechanics, the authors of \cite[Chapter 10]{LZ21} developed a refined mathematical model--namely, equation \eqref{1.3} with $f=g=h=0$--to  describe wave phenomena disturbed by random noise.

Controllability problems of deterministic hyperbolic equations have been extensively studied for several decades. It is impossible to provide a comprehensive list of references of this well-established theory. We refer readers to \cite{L88-1,L88-2,Z10,Z06} and the rich references therein for this topic. In contrast, the controllability theory for their stochastic counterparts remains significantly less developed.

To the best of our knowledge, \cite{LZ19} is the first work concerning the exact controllability of stochastic hyperbolic equations. In In that paper, the authors established boundary exact controllability for a refined stochastic hyperbolic equation under the key assumption that the control in the diffusion term of the first equation in \eqref{1.3} does not affect the drift term (i.e., $a_{5}=0$). This restrictive assumption was later removed in \cite{LL24}, where the authors proved exact controllability for refined stochastic hyperbolic equations with $(a^{jk}(\cdot))_{1\leq j,k\leq n}=I_{n}$. 
However, both \cite{LZ19} and \cite{LL24} require the waiting time $T$ (the time needed for the system to reach the desired target state) to be substantially longer than that required for deterministic wave equations. 

In the present work, we make two progresses on the study of controllability properties of the refined stochastic hyperbolic equations. First, we establish  internal exact controllability for system \eqref{1.3}. Secondly, we demonstrate that the required waiting time coincides exactly with that of the internal exact controllability for deterministic wave equations (cf. \cite{FLL23}), eliminating the need for additional waiting time typically required in stochastic settings.

\ms

To present our main result, we first introduce the following assumption  on the coefficients $(a^{jk}(\cdot))_{n\t n}$: 
\begin{condition}\label{con1.1}
There exists a positive function $\vp(\cdot)\in C^{2}(\ol{G})$ satisfying that
\begin{enumerate}
\item[(i)] For some constant $s_{0}\geq0$, and
$\forall (x,\xi)\in\ol{G}\t \dbR^{n}$, it holds that\vspace{-2mm}
\begin{equation}\label{1.8}
\sum_{j,k=1}^{n}\Big\{\sum_{j',k'=1}^{n}\Big[2a^{jk'}(a^{j'k}\vp_{x_{j'}})_{x_{k'}}
-a_{x_{k'}}^{jk}a^{j'k'}\vp_{x_{j'}}\Big]\Big\}\xi_{j}\xi_{k}\geq s_{0}\sum_{j,k=1}^{n}
a^{jk}\xi_{j}\xi_{k};
\end{equation}
\item[(ii)] There is no critical points of $\vp(\cdot)$ in $\ol{G}$, i.e.,\vspace{-2mm}
\begin{equation}\label{1.8-2}
\min_{x\in\ol{G}}|\nabla\vp(x)|>0.
\end{equation}
\end{enumerate}
\end{condition}
\begin{remark}
Condition \ref{con1.1} serves to ensure the existence of an appropriate weight function for deriving the required Carleman estimate for backward stochastic hyperbolic equations.  It is a kind of  pseudoconvex condition (see \cite[XXVIII]{Hormander09}).  Specifically, if the matrix $(a^{jk}(\cdot))_{1\leq j,k\leq n}$ equals the $n\times n$ identity matrix $I_n$, then for some $x_0 \in \mathbb{R}^n \setminus \overline{G}$, one can choose $\varphi(x) = |x - x_0|^2$ to satisfy Condition \ref{con1.1} with $s_0 = 4$.
\end{remark}

For a function $\vp(\cdot)$ satisfying Condition \ref{con1.1}, we define the subset \vspace{-2mm}
\begin{equation}\label{1.9}
\Gamma_{0}\=\Big\{x\in\Gamma~\Big|~\sum_{j,k=1}^{n}a^{jk}(x)\vp_{x_{j}}(x)\nu_{k}(x)>0\Big\},\vspace{-2mm}
\end{equation}
where $\nu(x)\=(\nu_{1}(x),\cdots,\nu_{n}(x))$ denotes the unit outward normal vector to $G$ at the point $x\in\Gamma$. We subsequently introduce the cylindrical domain $\Si_0\= (0, T)\t \G_0$.

One can readily verify that for any function $\vp(\cdot)\in C^{2}(\ol{G})$ satisfying \eqref{1.8}, the transformed function $\hat{\vp}(\cdot)=a\vp(\cdot)+b$ maintains this property for all $a\geq1$ and $b\in\dbR$, with $s_{0}$ replaced by $as_{0}$. Crucially, such affine transformations preserve the set $\Gamma_{0}$ defined in \eqref{1.9}. Consequently, without loss of generality, we may impose the following strengthened conditions throughout our analysis: \vspace{-2mm}
\begin{equation}\label{1.9-1}
\begin{cases}
\text{\eqref{1.8} holds with $s_{0}\geq4$},\\
\frac{1}{4}\sum_{j,k=1}^{n}a^{jk}(x)\vp_{x_{j}}(x)\vp_{x_{k}}(x)\geq
\vp(x)>0,\q \forall x\in\ol{G}.
\end{cases}
\end{equation}

For any given $\delta>0$, let\vspace{-2mm}
\begin{equation}\label{1.10}
\cO_{\delta}(\Gamma_{0})\=\{x\in \dbR^{n}~|~\dist(x,\Gamma_{0})<\delta \},
\end{equation}
and\vspace{-2mm}
\begin{equation}\label{1.11}
G_{0}\=\cO_{\delta}(\Gamma_{0})\cap G.
\end{equation}
Put\vspace{-2mm}
\begin{equation}\label{1.12}
M_{0}\=\min_{x\in\ol{G}}\sum_{j,k=1}^{n}a^{jk}\vp_{x_{j}}\vp_{x_{k}},
\q M_{1}\=\max_{x\in\ol{G}}\sum_{j,k=1}^{n}a^{jk}\vp_{x_{j}}\vp_{x_{k}} \vspace{-2mm}
\end{equation}
and\vspace{-2mm}
\begin{equation}\label{1.14}
T_{0}\=\sqrt{M_{1}}.
\end{equation}

\begin{remark}
For the special case where the coefficient matrix reduces to the identity matrix, i.e., \linebreak $(a^{jk}(\cdot))_{1\leq j,k\leq n} =I_{n}$, we can choose $\vp(x)=|x-x_{0}|^{2}$ with $x_{0}\in\dbR^{n}\setminus\ol{G}$. Then \vspace{-2mm}
$$
T_{0}=2\max_{x\in\ol{G}}|x-x_{0}|.
$$
\end{remark}

Now we give the main  result of this paper.
\begin{theorem}\label{thm1.1}
Under Condition \ref{con1.1}, system \eqref{1.3} is exactly controllable for any time  $T>T_{0}$. More precisely, given any initial states  $(y_{0},\hat{y}_{0})\in H^{1}_{0}(G)\t L^{2}(G)$ and target states $(y_{1},\hat{y}_{1})\in L_{\cF_{T}}^{2}(\Omega;H^{1}_{0}(G))\t L_{\cF_{T}}^{2}(\Omega;L^{2}(G))$, there exist controls $f\in L_{\dbF}^{2}(0,T;H^{1}_{0}(G))$, $g\in L_{\dbF}^{2}(0,T;L^{2}(G))$, and $h\in L_{\dbF}^{2}(0,T;L^{2}(G))$ such that the corresponding solution $(y,\hat{y})$ to \eqref{1.3} satisfies that $(y(T),\hat{y}(T))=(y_{1},\hat{y}_{1})$.
\end{theorem}

To establish Theorem \ref{thm1.1}, we employ a standard duality argument (cf. \cite[Section 7.3]{LZ21}), which reduces the proof to deriving an observability estimate for the following backward stochastic hyperbolic equation:\vspace{-2mm}
\begin{equation}\label{1.4}
\begin{cases}
dz=\hat{z}dt+ZdW(t) & \text{in}\quad Q,\\
d\hat{z}-\sum_{j,k=1}^{n}(a^{jk}z_{x_{j}})_{x_{k}}dt=(a_{1}z+a_{2}Z-a_{3}\h{Z})dt+\h{Z}dW(t) &\text{in}\quad Q,\\
z=0 & \text{on} \quad \Sigma,\\
z(T)=z^{T},\quad \hat{z}(T)=\hat{z}^{T} & \text{in}\quad G,
\end{cases}\vspace{-2mm}
\end{equation}
where   $(z^{T},\hat{z}^{T})\in L^{2}_{\cF_{T}}(\Omega;L^{2}(G)\times H^{-1}(G))$.

By the classical well-posedness result for backward stochastic evolution equations (see, e.g., \cite[Theorem 4.10]{LZ21}), we conclude that for any terminal data $(z^{T},\hat{z}^{T})\in L^{2}_{\cF_{T}}(\Omega;L^{2}(G)\times H^{-1}(G))$, the system \eqref{1.4} admits a unique weak solution
$(z,\hat{z},Z,\h{Z})$ in
$L^{2}_{\dbF}(\Omega;C([0,T];L^{2}(G)))\times
L^{2}_{\dbF}(\Omega;C([0,T];H^{-1}(G)))\times
L^{2}_{\dbF}(0,T;L^{2}(G))\times
L^{2}_{\dbF}(0,T;H^{-1}(G))$. 

According to \cite[Theorem 7.16]{LZ21}, Theorem \ref{thm1.1} is equivalent to the following  internal observability estimate.

\begin{theorem}\label{thm1.2}
Assume that Condition \ref{con1.1} holds. Then for any terminal time   $T>T_{0}$, there exists a positive constant $C >0$ such that for all terminal data   $(z^{T},\hat{z}^{T})\in L^{2}_{\cF_{T}}(\Omega;L^{2}(G)\times H^{-1}(G))$ and the corresponding solution to system \eqref{1.4}, the following observability inequality holds: \vspace{-2mm}
\begin{equation}\label{1.15}
\begin{aligned}
&\|(z^{T},\hat{z}^{T})\|^{2}_{L^{2}_{\cF_{T}}(\Omega;L^{2}(G)\times H^{-1}(G))}\\
& \leq  C \Big( \dbE\int_{0}^{T}\int_{G_{0}} z^{2}dxdt
+\|a_{4}z+Z\|^{2}_{L^{2}_{\dbF}(0,T;L^{2}(G))}
+\|a_{5}\hat{z}+\h{Z}\|^{2}_{L^{2}_{\dbF}(0,T;H^{-1}(G))}\Big).
\end{aligned}\vspace{-2mm}
\end{equation}
\end{theorem}

Hereafter, unless otherwise specified, we denote by $C=C(\vp, T, (a^{jk})_{n\t n}, G, G_0,  a_1, a_2, a_3,a_4,a_5)$ a generic positive constant that may vary from line to line.

In this paper, we prove Theorem \ref{thm1.2} by establishing an internal Carleman estimate for the  backward stochastic hyperbolic equation \eqref{1.4}.  
A usual way to establish the internal observability estimate for hyperbolic equation  is to combine a suitable boundary Carleman estimate, the multiplier method and a suitable energy estimate (e.g., \cite{FYZ07,FLL23,Imanuvilov02}). This strategy is generalized to the stochastic setting in \cite{FLLZ16}, in which the authors  establish an internal Carleman estimate for stochastic hyperbolic equations.

Following this general framework, one might attempt to establish an internal observability estimate for our backward stochastic hyperbolic equation \eqref{1.4} by utilizing the boundary Carleman estimates developed in \cite{LL24,LZ21} for backward stochastic hyperbolic equations. However, as previously noted, the waiting time required in \cite{LL24,LZ21} significantly exceeds the time $T_0$ defined in \eqref{1.14}. This discrepancy prevents us from directly adapting the proof techniques of \cite{FLLZ16} to prove Theorem \ref{thm1.2}.

To circumvent this obstacle, we develop an enhanced duality argument and derive the desired internal Carleman estimate through an alternative approach based on boundary Carleman estimates for deterministic hyperbolic operators and suitable energy estimate for backward stochastic hyperbolic equations.

\ss

The rest of this paper is organized as follows. In Section \ref{sec2}, we establish a Carleman estimate for stochastic hyperbolic equations in $L^2$-norm based on a Carleman estimate for deterministic wave operators in $H^1$-norm. In Section \ref{sec3}, we prove three energy estimates for the solution to system \eqref{1.4}, and then prove our main result.  

\section{An $L^2$-Carleman estimate for a hyperbolic operator}\label{sec2}

In this section, we establish an $L^{2}$-Carleman estimate for the following hyperbolic operator:\vspace{-2mm}
$$
\cA\=\frac{\partial^{2}}{\partial t^{2}}-\sum_{j,k=1}^{n}\frac{\partial}{\partial{x_{k}}}\Big(
a^{jk}\frac{\partial}{\partial{x_{j}}}\Big).\vspace{-2mm}
$$

\subsection{An auxiliary optimal control problem}\label{ssec2.1}

In this subsection,   borrowing some idea from  \cite{FLLZ16}, we introduce an auxiliary optimal control problem that plays a crucial role in establishing the Carleman estimate for equation \eqref{1.4}. 

\ss

For any $K>1$, we choose $\rho^{K}\equiv\rho^{K}(x)\in C^{2}(\ol{G})$ such
that $\min_{x\in\ol{G}}\rho^{K}(x)=1$ and\vspace{-2mm}
\begin{equation}\label{2.3}
\rho^{K}(x)=
\begin{cases}
1,& x\in G_{0},\\
K, & \dist(x,G_{0})\geq\frac{1}{\ln K}.
\end{cases}\vspace{-2mm}
\end{equation}

For some given $c_{1}\in\Big(\frac{\sqrt{M_{1}}}{T},1\Big)$ and parameters $\l,\mu>0$, we define\vspace{-2mm}
\begin{equation}\label{2.1}
\th(t,x)=e^{\l\phi(t,x)},
\q \phi(t,x)=e^{\mu\si(t,x)},
\q \si(t,x)=\vp(x)-c_{1}\Big(t-\frac{T}{2}\Big)^{2}.
\end{equation}

Let $v\in L_{\dbF}^{2}(\Omega;C([0,T];L^{2}(G)))$ satisfy $v(0)=v(T)=0$, $v|_{\Si}=0$ and $\supp v\subset[0,T]\t(\ol{G}\setminus G_{0})$. For any integer $m\geq3$, let
$\D t=\frac{T}{m}$, and set\vspace{-2mm}
\begin{equation}\label{2.4}
v_{m}^{j}=v_{m}^{j}(x)=v(j\D t,x),
\q \phi_{m}^{j}=\phi_{m}^{j}(x)=\phi(j\D t,x),
\q 0\leq j \leq m,
\end{equation}
where $\phi$ is given in \eqref{2.1}. Consider the following system:\vspace{-2mm}
\begin{equation}\label{2.5}
\begin{cases}
\dbE\Big(\frac{w_{m}^{j+1}-2w_{m}^{j}+w_{m}^{j-1}}{\D t^2}\Big|\cF_{j\D t}\Big)
-\sum_{j_{1},j_{2}=1}^{n}\partial_{x_{j_{2}}}(a^{j_{1}j_{2}}\partial_{x_{j_{2}}}w_{m}^{j})\\
=\dbE\Big(\frac{r_{1m}^{j+1}-r_{1m}^{j}}{\D t}\Big|\cF_{j\D t}\Big)+\lambda v_{m}^{j}e^{2\l\phi_{m}^{j}}+r_{m}^{j},
\q 1\leq j\leq m-1 &\text{in}\q G,\\
w_{m}^{j}=0, \q 0\leq j\leq m &\text{on} \q \Gamma,\\
w_{m}^{0}=w_{m}^{m}=r_{m}^{0}=r_{m}^{m}=0, \q r_{1m}^{0}=r_{1m}^{1} & \text{in}\q G.
\end{cases}
\end{equation}
Here $(r_{1m}^{j},r_{m}^{j})\in (L_{\cF_{j\D t}}^{2}(\Omega;L^{2}(G)))^{2}$
are controls with $0\leq j\leq m$. The set of admissible sequences for \eqref{2.5} is defined by\vspace{-2mm}
\begin{equation}\label{2.6}
\begin{aligned}
\sA_{ad}\=
&\Big\{
\{(w_{m}^{j},r_{1m}^{j},r_{m}^{j})\}_{j=0}^{m}~\Big|~
(w_{m}^{j},r_{1m}^{j},r_{m}^{j})\in L_{\cF_{j\D t}}^{2}(\Omega;H^{1}_{0}(G))\t (L_{\cF_{j\D t}}^{2}(\Omega;L^{2}(G)))^{2}\\
&\q\text{fulfills \eqref{2.5} and}
~\supp (w_{m}^{j},r_{1m}^{j},r_{m}^{j})\subset (\ol{G}\setminus G_{0})^{3}
~\text{for}~0\leq j\leq m
\Big\}.
\end{aligned}
\end{equation}
Since $\{(0,0,-\lambda v_{m}^{j}e^{2\lambda\phi_{m}^{j}})\}_{j=0}^{m}\in \sA_{ad}$, we know
that $\sA_{ad}\neq \emptyset$.

Next we define a cost functional\vspace{-2mm}
\begin{eqnarray}\label{2.7}  \cJ\big(\{(w_{m}^{j},r_{1m}^{j},r_{m}^{j})\}_{j=0}^{m}\big)&=&\frac{\D t}{2}\dbE\sum_{j=1}^{m-1}\Big[\int_{G}|w_{m}^{j}|^{2}e^{-2\l\phi_{m}^{j}}dx
+\int_{G}\rho^{K}\frac{|r_{1m}^{j}|^{2}}{\l^{2}}
e^{-2\l\phi_{m}^{j}}dx \nonumber\\
&&\qq\qq\q +K\int_{G}|r_{m}^{j}|^{2}dx\Big]+\frac{\D t}{2}\dbE\int_{G}\rho^{K}\frac{|r_{1m}^{m}|^{2}}{\lambda^{2}}e^{-2\l \phi_{m}^{m}}dx,
\end{eqnarray}
and consider the following optimal control problem: 

Find a
$\{(\hat{w}_{m}^{j},\hat{r}_{1m}^{j},\hat{r}_{m}^{j})\}_{j=0}^{m}\in\sA_{ad}$ such that
\begin{equation}\label{2.8}
\cJ\big(\{(\hat{w}_{m}^{j},\hat{r}_{1m}^{j},\hat{r}_{m}^{j})\}_{j=0}^{m}\big)
=\min_{\{(w_{m}^{j},r_{1m}^{j},r_{m}^{j})\}_{j=0}^{m}\in\sA_{ad}}\cJ\big(\{(w_{m}^{j},r_{1m}^{j},r_{m}^{j})\}_{j=0}^{m}\big).
\end{equation}
For any $\{(w_{m}^{j},r_{1m}^{j},r_{m}^{j})\}_{j=0}^{m}\in\sA_{ad}$, standard elliptic regularity theory yields $w_{m}^{j}\in L_{\cF_{j\D t}}^{2}(\Omega;H^{2}(G)\cap H^{1}_{0}(G))$
for $0\leq j\leq m$. Furthermore, we have the following result.
\begin{proposition}\label{prop2.1}
For any $m\geq3$ and $K>1$,  the optimal control problem \eqref{2.8} admits a unique solution
$\{(\hat{w}_{m}^{j},\hat{r}_{1m}^{j},\hat{r}_{m}^{j})\}_{j=0}^{m}\in\sA_{ad}$, which depends on 
$K$. Defining the adjoint state  \vspace{-2mm}
\begin{equation}\label{2.9}
p_{m}^{j}=p_{m}^{j}(x)\=K\hat{r}_{m}^{j}(x), \q 0\leq j\leq m,
\end{equation}
we have \vspace{-2mm}
\begin{equation}\label{2.10}
\begin{cases}
\ds\hat{w}_{m}^0=\hat{w}_{m}^{m}=p_{m}^{0}=p_{m}^{m}=0\q\text{in}\q G,\\
\ns\ds \hat{w}_{m}^{j},p_{m}^{j}\in L^{2}_{\cF_{j\D t}}(\Omega;H^{2}(G)\cap H^{1}_{0}(G)),\q 1\leq j\leq m-1.
\end{cases}
\end{equation}
The solution satisfies the following optimality conditions:\vspace{-2mm}
\begin{equation}\label{2.11}
\frac{p_{m}^{j}-p_{m}^{j-1}}{\D t}+\rho^{K}\frac{\hat{r}_{1m}^{j}}{\l^{2}}e^{-2\l\phi_{m}^{j}}=0
\q\text{in}\q G,\q 1\leq j\leq m,
\end{equation}
and\vspace{-2mm}
\begin{equation}\label{2.12}
\begin{cases}
\dbE\Big(\frac{p_{m}^{j+1}-2p_{m}^{j}+p_{m}^{j-1}}{\D t^{2}}\Big|\cF_{j\D t}\Big)
-\sum_{j_{1},j_{2}=1}^{n}\partial_{x_{j_{2}}}(a^{j_{1}j_{2}}\partial_{x_{j_{1}}}p_{m}^{j})
+e^{-2\l\phi_{m}^{j}}\hat{w}_{m}^{j}=0 & \text{in}\q G,\\
p_{m}^{j}=0 & \text{on}\q \Gamma.
\end{cases}
\end{equation}
Moreover, there exists a positive constant 
$C=C(\l,K)$, independent of $m$, such that the following estimates hold: 
\begin{equation}\label{2.13}
\D t\dbE\sum_{j=1}^{m-1}\int_{G}\big(|\hat{w}_{m}^{j}|^{2}+|\hat{r}_{1m}^{j}|^{2}
+K|\hat{r}_{m}^{j}|^{2}\big)dx+\D t\dbE\int_{G}|\hat{r}_{1m}^{m}|^{2}dx\leq C,
\end{equation}
\begin{equation}\label{2.14}
\D t\dbE\sum_{j=0}^{m-1}\int_{G}\Big[\frac{(\hat{w}_{m}^{j+1}-\hat{w}_{m}^{j})^{2}}{\D t^{2}}
+\frac{(\dbE(\hat{r}_{1m}^{j+1}-\hat{r}_{1m}^{j}|\cF_{j\D t}))^{2}}{\D t^{2}}
+K\frac{(\hat{r}_{m}^{j+1}-\hat{r}_{m}^{j})^2}{\D t^{2}}\Big]dx\leq C,
\end{equation}
and\vspace{-2mm}
\begin{equation}\label{2.15}
\D t\dbE\sum_{j=1}^{m-1}\int_{G}\big(|\nabla\hat{w}_{m}^{j}|^{2}+|\nabla \hat{r}_{m}^{j}|^{2}\big)dx\leq C.
\end{equation}
\end{proposition}

Before  proving Proposition \ref{prop2.1}, we  recall the following known result.
\begin{lemma}{\rm\cite[Proposition 3.5]{FYZ07}}\label{lem-a1}
For any $m\geq3$, and $q_{m}^{j},w_{m}^{j}\in\dbC$, $0\leq j\leq m$ satisfying $q_{0}^{m}=q_{m}^{m}=0$, one has
\begin{equation}\label{a1}
\begin{aligned}
-\sum_{j=1}^{m-1}q_m^j\frac{w_m^{j+1}-2w_m^j+w_m^{j-1}}{\tau^2}   =\sum_{j=0}^{m-1}\frac{q_m^{j+1}-q_m^j}{\tau}\frac{w_m^{j+1}-w_m^j}{\tau} =\sum_{j=1}^m\frac{q_m^j-q_m^{j-1}}{\tau}\frac{w_m^j-w_m^{j-1}}{\tau}.
\end{aligned}
\end{equation}
\end{lemma}

{\it Proof of Proposition \ref{prop2.1}}. 
The argument is lengthy and will be established through five steps.

\ss

{\bf Step 1.} In this step, we show that the above optimal control problem has a unique solution. 

It is an easy matter to see that the functional $\cJ(\cdot)$ is weakly lower semi-continuous, strictly convex, and coercive in $\prod_{j=0}^{m} (L^{2}_{\cF_{j\D t}}(\Omega;L^{2}(G)))^{3}$. The coercivity property of $\cJ(\cdot)$ guarantees the existence of a bounded minimizing sequence  $
\big\{\{(w_{m}^{j,k},r_{1m}^{j,k},r_{m}^{j,k})\}_{j=0}^{m}\big\}_{k=1}^{\i}
\subset\sA_{ad}$. 
By weak compactness, there exists a subsequence of $\big\{\{(w_{m}^{j,k},r_{1m}^{j,k},r_{m}^{j,k})\}_{j=0}^{m}\big\}_{k=1}^{\i}$ that converges weakly to some limit 
$$
\{(\hat{w}_{m}^{j},\hat{r}_{1m}^{j},\hat{r}_{m}^{j})\}_{j=0}^{m}\in\sA_{ad}
\q\text{in}\q
\prod_{j=0}^{m}
(L_{\cF_{j\D t}}^{2}(\Omega;L^{2}(G)))^{3}.
$$
Elliptic regularity theory further yields the improved regularity:
$$
\{(\hat{w}_{m}^{j},\hat{r}_{1m}^{j},\hat{r}_{m}^{j})\}_{j=0}^{m}\in
\prod_{j=0}^{m}\Big[L_{\cF_{j\D t}}^{2}(\Omega;H^{1}_{0}(G))\t
(L_{\cF_{j\D t}}^{2}(\Omega;L^{2}(G)))^{2}\Big].
$$
An application of Mazur's Lemma shows that $(\hat{w}_{m}^{j},\hat{r}_{1m}^{j},\hat{r}_{m}^{j})$
can be obtained as the strong limit of convex combinations from the weakly convergent subsequence. This implies that
$\supp(\hat{w}_{m}^{j},\hat{r}_{1m}^{j},\hat{r}_{m}^{j})\subset(\ol{G}\setminus G_{0})^{3}$ for $0\leq j\leq m$. The weak lower semi-continuity of $\cJ(\cdot)$   ensures that $(\hat{w}_{m}^{j},\hat{r}_{1m}^{j},\hat{r}_{m}^{j})$   indeed solves the optimal control problem \eqref{2.8}, while the strict convexity guarantees uniqueness. Finally, equation \eqref{2.9} combined with the definition of 
$\sA_{ad}$ yields $\hat{w}_{m}^0=\hat{w}_{m}^{m}=p_{m}^{0}=p_{m}^{m}=0$ in $G$.

\ss

{\bf Step 2.} For each $0\leq j\leq m$, fix $w_{0m}^{j}\in L^{2}_{\cF_{j\D t}}(\Omega;H^{2}(G)\cap H^{1}_{0}(G))$ and $w_{1m}^{j}\in L^{2}_{\cF_{j\D t}}(\Omega;L^{2}(G))$ with  the conditions \vspace{-2mm}
$$w_{0m}^{0}=w_{0m}^{m}=0,\qq w_{1m}^{0}=w_{1m}^{1}$$ 
and\vspace{-2mm}
$$\supp(w_{0m}^{j},w_{1m}^{j})\subset(\ol{G}\setminus G_{0})^{2}.$$ 
For any $(\gamma_{0},\gamma_{1})\in\dbR^{2}$, the triplet
$\{(\hat{w}_{m}^{j}+\g_{0}w_{0m}^{j},\hat{r}_{1m}^{j}+\g_{1}w_{1m}^{j},r_{m}^{j})\}
_{j=0}^{m}$ belongs to $\sA_{ad}$, 
where
\begin{equation}\label{a2}
\begin{cases}
r_{m}^{j}=\dbE\Big(\frac{\hat{w}_{m}^{j+1}-2\hat{w}_{m}^{j}+\hat{w}_{m}^{j-1}}{\D t^{2}}\Big|\cF_{j\D t}\Big)
-\sum_{j_{1},j_{2}=1}^{n}\partial_{x_{j_{2}}}(a^{j_{1}j_{2}}\partial_{x_{j_{1}}}\hat{w}_{m}^{j})
-\mathbb{E}\Big(\frac{\hat{r}_{1m}^{j+1}-\hat{r}_{1m}^{j}}{\D t}\Big|\mathcal{F}_{j\D t}\Big)\\
\qq\q-\lambda v_{m}^{j}e^{2\lambda\phi_{m}^{j}}+\g_{0}\Big[
\dbE\Big(\frac{w_{0m}^{j+1}-2w_{0m}^{j}+w_{0m}^{j-1}}{\D t^{2}}\Big|\cF_{j\D t}\Big)
-\sum_{j_{1},j_{2}=1}^{n}\partial_{x_{j_{2}}}(a^{j_{1}j_{2}}\partial_{x_{j_{1}}}w_{0m}^{j})\Big] \\ \qq\q-\g_{1}\dbE\Big(\frac{w_{1m}^{j+1}-w_{1m}^{j}}{\D t}\Big|\cF_{j\D t}\Big),\quad
1\leq j\leq m-1; \\
r_{m}^{0}=r_{m}^{m}=0.
\end{cases}
\end{equation}
Define the function $\cH:\dbR^{2}\to\dbR$ by
\begin{equation}\label{a3}
\cH(\g_{0},\g_{1})=\cJ\big(\{(\hat{w}_{m}^{j}+\g_{0}w_{0m}^{j},\hat{r}_{1m}^{j}+\g_{1}w_{1m}^{j},r_{m}^{j})\}
_{j=0}^{m}\big).
\end{equation}
Since $\cH$ attains its minimum at  $(0,0)$, it follows that $\nabla \cH(0,0)=0$. 

From $\frac{\partial \cH(0,0)}{\partial \g_{1}}=0$, we deduce
\begin{equation}\label{a4}
\mathbb{E}\sum_{j=1}^{m-1}\int_{G}\Big(K\frac{\hat{r}_{m}^{j}-\hat{r}_{m}^{j-1}}{\D t}
+\rho^{K}\frac{\hat{r}_{1m}^{j}}{\l^{2}}e^{-2\lambda\phi_{m}^{j}}\Big)w_{1m}^{j}dx=0,
\end{equation}
which yields \eqref{2.11}. Similarly,    $\frac{\partial \cH(0,0)}{\partial \g_{0}}=0$ implies
\begin{equation}\label{a5}
\begin{aligned}
&\mathbb{E}\sum_{j=1}^{m-1}\int_{G}\bigg\{K\hat{r}_{m}^{j}\Big[\mathbb{E}\Big(
\frac{w_{0m}^{j+1}-2w_{0m}^{j}+w_{0m}^{j-1}}{\D t^{2}}\Big|\mathcal{F}_{j\D t}\Big)
-\sum_{j_{1},j_{2}=1}^{n}\partial_{x_{j_{2}}}(a^{j_{1}j_{2}}\partial_{x_{j_{1}}}w_{0m}^{j})\Big]\\
&\qq\qq\q+\hat{w}_{m}^{j}w_{0m}^{j}e^{-2\lambda\phi_{m}^{j}}\bigg\}dx=0.
\end{aligned}
\end{equation}
Combined with the  conditions $p_{m}^{0}=p_{m}^{m}=w_{0m}^{0}=w_{0m}^{m}=0$ in $G$, this gives \eqref{2.12}. By the regularity theory for the solutions to elliptic equation, we conclude that $\hat{w}_{m}^{j},p_{m}^{j}\in L^{2}_{\cF_{j\D t}}(\Omega;H^{2}(G)\cap H^{1}_{0}(G))$, $1\leq j\leq m-1$.

\ss

{\bf Step 3.} In this step, we establish the estimate \eqref{2.13}. Substituting $\{w_{0m}^{j}\}_{j=0}^{m}$ with $\{\hat{w}_{m}^{j}\}_{j=0}^{m}$ in \eqref{a5} and recalling that $p_{m}^{j}=K\hat{r}_{m}^{j}$, we  derive
\begin{eqnarray}\label{a6}
0& =&\mathbb{E}\sum_{j=1}^{m-1}\int_{G}\Big[\mathbb{E}\Big(\frac{\hat{w}_{m}^{j+1}-2\hat{w}_{m}^{j}+\hat{w}_{m}^{j-1}}{\D t^{2}} \Big| \mathcal{F}_{j\D t}\Big)-\sum_{j_{1},j_{2}=1}^{n}\partial_{x_{j_{2}}}(a^{j_{1}j_{2}}\partial_{x_{j_{1}}}\hat{w}_{m}^{j})\nonumber\\
&& \qq\qq\qq-\mathbb{E}\Big(\frac{\hat{r}_{1m}^{j+1}-\hat{r}_{1m}^{j}}{\D t} \Big| \mathcal{F}_{j\D t}\Big)-\lambda v_{m}^{j}e^{2\lambda\phi_{m}^{j}}-\hat{r}_{m}^{j} \Big]p_{m}^{j}dx \nonumber\\
& =&\mathbb{E}\sum_{j=1}^{m-1}\int_{G}\Big[\mathbb{E}\Big(\frac{p_{m}^{j+1}-2p_{m}^{j}+p_{m}^{j-1}}{\D t^{2}} \Big| \mathcal{F}_{j\D t}\Big)-\sum_{j_{1},j_{2}=1}^{n}\partial_{x_{j_{2}}}(a^{j_{1}j_{2}}\partial_{x_{j_{1}}}p_{m}^{j})\Big]\hat{w}_{m}^{j}dx\\
&& \q +\mathbb{E}\sum_{j=1}^{m}\int_{G}\frac{p_{m}^{j}-p_{m}^{j-1}}{\D t}\hat{r}_{1m}^{j}dx
-\mathbb{E}\sum_{j=1}^{m-1}\int_{G}\Big(\lambda v_{m}^{j}e^{2\lambda\phi_{m}^{j}}+\hat{r}_{m}^{j}\Big)p_{m}^{j}dx \nonumber\\
& =&-\mathbb{E}\sum_{j=1}^{m-1}\Big[\int_{G}|\hat{w}_{m}^{j}|^{2}e^{-2\lambda\phi_{m}^{j}}dx
+\int_{G}\rho^{K}\frac{|\hat{r}_{1m}^{j}|^{2}}{\lambda^{2}}e^{-2\lambda\phi_{m}^{j}}dx
+K\int_{G}|\hat{r}_{m}^{j}|^{2}dx\Big]\nonumber\\
&& \q-\mathbb{E}\int_G\rho^{K}\frac{|\hat{r}_{1m}^m|^2}{\lambda^2}e^{-2\lambda\phi_m^m}dx
-\lambda\mathbb{E}\sum_{j=1}^{m-1}\int_{G}v_{m}^{j}e^{2\lambda\phi_{m}^{j}}p_{m}^{j}dx. \nonumber
\end{eqnarray}
Combining \eqref{2.11} and \eqref{a6}, there exists a constant 
$C=C(K,\l)>0$,  independent of $m$, such that
\begin{equation}\label{a7}
\begin{aligned}
&\mathbb{E}\sum_{j=1}^{m-1}\Big[\int_{G}|\hat{w}_{m}^{j}|^{2}e^{-2\lambda\phi_{m}^{j}}dx
+\int_{G}\rho^{K}\frac{|\hat{r}_{1m}^{j}|^{2}}{\lambda^{2}}e^{-2\lambda\phi_{m}^{j}}dx
+K\int_{G}|\hat{r}_{m}^{j}|^{2}dx\Big]\\
&+\mathbb{E}\int_{G}\rho^{K}\frac{|\hat{r}_{1m}^{m}|^{2}}{\lambda^{2}}e^{-2\lambda\phi_{m}^{m}}dx
\leq C\mathbb{E}\sum_{j=1}^{m-1}\int_{G}|v_{m}^{j}|^{2}e^{2\lambda\phi_{m}^{j}}dx,
\end{aligned}
\end{equation}
which concludes the proof of \eqref{2.13}.

\ss

{\bf Step 4.} In this step we prove \eqref{2.14}.
We begin by recalling equation  \eqref{2.5}, which yields
\begin{equation}\label{a8}
\begin{aligned}
0&=\dbE\sum_{j=1}^{m-1}\int_{G}\bigg[\dbE\Big(\frac{\hat{w}_{m}^{j+1}-2\hat{w}_{m}^{j}
+\hat{w}_{m}^{j-1}}{\D t^{2}}\Big|\cF_{j\D t}\Big)-\sum_{j_{1},j_{2}=1}^{n}\partial_{x_{j_{2}}}
(a^{j_{1}j_{2}}\partial_{x_{j_{1}}}\hat{w}_{m}^{j})\\
&\qq\qq\qq-\dbE\Big(\frac{\hat{r}_{1m}^{j+1}-\hat{r}_{1m}^{j}}{\D t}\Big|\cF_{j\D t}\Big)-\l v_{m}^{j}e^{2\l\phi_{m}^{j}}-\hat{r}_{m}^{j}\bigg]\dbE\Big(\frac{p_{m}^{j+1}-2p_{m}^{j}
+p_{m}^{j-1}}{\D t^{2}}\Big|\cF_{j\D t}\Big)dx.
\end{aligned}
\end{equation}
From equation \eqref{2.12} and the conditions $p_{m}^{0}=p_{m}^{m}=\hat{w}_{m}^{0}=\hat{w}_{m}^{m}=0$ in $G$, we  obtain that\vspace{-2mm}
\begin{equation}\label{a9}
\begin{aligned}
&\mathbb{E}\sum_{j=1}^{m-1}\int_{G}\mathbb{E}\Big({\frac{\hat{w}_{m}^{j+1}-2\hat{w}_{m}^{j}+\hat{w}_{m}^{j-1}}{\D t^{2}}} \Big|\cF_{j\D t}\Big)\mathbb{E}\Big({\frac{p_{m}^{j+1}-2p_{m}^{j}+p_{m}^{j-1}}{\D t^{2}}} \Big|\cF_{j\D t}\Big)dx\\
&=\mathbb{E}\sum_{j=1}^{m-1}\int_{G}\hat{w}_{m}^{j}\bigg\{\sum_{j_{1},j_{2}=1}^{n}\partial_{x_{j_{2}}}
\Big[a^{j_{1}j_{2}}\partial_{x_{j_{1}}}\mathbb{E}\Big(\frac{p_{m}^{j+1}-2p_{m}^{j}+p_{m}^{j-1}}{\D t^{2}}
\Big|\mathcal{F}_{j\D t}\Big)\Big] \\
&\qq\qq\qq\qq-\frac{\mathbb{E}(\hat{w}_{m}^{j+1}|\mathcal{F}_{j\D t})e^{-2\lambda\phi_{m}^{j+1}} -2\hat{w}_{m}^{j}e^{-2\lambda\phi_{m}^{j}}+\hat{w}_{m}^{j-1}e^{-2\lambda\phi_{m}^{j-1}}}{\D t^{2}}\bigg\}dx.
\end{aligned}
\end{equation}
Since $\hat{w}_m^j|_{\Gamma}=p_{m}^{j}|_{\Gamma}=0$ for $0\leq j\leq m$, integration by parts gives\vspace{-2mm}
\begin{equation}\label{a10}
\begin{aligned}
&\mathbb{E}\sum_{j=1}^{m-1}\int_{G}\sum_{j_{1},j_{2}=1}^{n}\partial_{x_{j_{2}}}
\big(a^{j_{1}j_{2}}\partial_{x_{j_{1}}}\hat{w}_{m}^{j}\big)
\mathbb{E}\Big(\frac{p_{m}^{j+1}-2p_{m}^{j}+p_{m}^{j-1}}{\D t^{2}}\Big|\mathcal{F}_{j\D t}\Big)dx\\
&=\mathbb{E}\sum_{j=1}^{m-1}\int_{G}\hat{w}_{m}^{j}\sum_{j_{1},j_{2}=1}^{n}\partial_{x_{j_{2}}}
\Big[a^{j_{1}j_{2}}\partial_{x_{j_{1}}}\mathbb{E}\Big(\frac{p_{m}^{j+1}-2p_{m}^{j}+p_{m}^{j-1}}{\D t^{2}}
\Big|\mathcal{F}_{j\D t}\Big)\Big]dx.
\end{aligned}
\end{equation}
Combining equations (\ref{a8})-(\ref{a10}), we obtain that\vspace{-2mm}
\begin{eqnarray}\label{a11} 
0 \3n& =&\3n-\mathbb{E}\sum_{j=1}^{m-1}\int_{G}\Big[\hat{w}_m^j
\frac{\dbE(\hat{w}_{m}^{j+1}|\mathcal{F}_{j\D t})e^{-2\lambda\phi_m^{j+1}}
-2\hat{w}_m^{j}e^{-2\lambda\phi_m^{j}}+\hat{w}_m^{j-1}e^{-2\lambda\phi_m^{j-1}}}{\D t^{2}}\\
&& \qq\qq\qq+\Big(\frac{\dbE(\hat{r}_{1m}^{j+1}|\cF_{j\D t})\!-\!\hat{r}_{1m}^{j}}{\D t}
\!+\!\lambda v_m^{j}e^{2\lambda\phi_m^{j}}\!+\!\hat{r}_m^j\Big)\mathbb{E}\Big(\frac{p_m^{j+1}\!-\!2p_m^j\!+\!p_m^{j-1}}{\D t^2}
\Big|\mathcal{F}_{j\D t}\Big)\Big]dx.\nonumber
\end{eqnarray}
Applying Lemma \ref{lem-a1} to the first term yields
\begin{equation}\label{a12}
\begin{aligned}
&-\mathbb{E}\sum_{j=1}^{m-1}\int_{G}\Big(\hat{w}_{m}^{j}\frac{\mathbb{E}(\hat{w}_{m}^{j+1}|\mathcal{F}_{j\D t})e^{-2\lambda\phi_{m}^{j+1}}
-2\hat{w}_{m}^{j}e^{-2\lambda\phi_{m}^{j}}+\hat{w}_{m}^{j-1}e^{-2\lambda\phi_{m}^{j-1}}}{\D t^{2}}\Big)dx\\
&=\mathbb{E}\sum_{j=0}^{m-1}\int_{G}\frac{(\hat{w}_{m}^{j+1}-\hat{w}_{m}^{j})}{\D t}\frac{(\hat{w}_{m}^{j+1}e^{-2\lambda\phi_{m}^{j+1}}
-\hat{w}_{m}^{j}e^{-2\lambda\phi_{m}^{j}})}{\D t}dx \\
&=\mathbb{E}\sum_{j=0}^{m-1}\int_{G}\Big[\frac{(\hat{w}_{m}^{j+1}-\hat{w}_{m}^{j})^{2}}{\D t^{2}}e^{-2\lambda\phi_{m}^{j}}
+\frac{(\hat{w}_{m}^{j+1}-\hat{w}_{m}^{j})}{\D t}\frac{(e^{-2\lambda\phi_{m}^{j+1}}-e^{-2\lambda\phi_{m}^{j}})}{\D t}\hat{w}_{m}^{j+1}\Big]dx.
\end{aligned}\vspace{-2mm}
\end{equation}
By Lemma \ref{lem-a1} and $p_{m}^{j}=K\hat{r}_{m}^{j}$, we obtain that
\begin{equation}\label{a13}
-\mathbb{E}\sum_{j=1}^{m-1}\int_{G}\Big[\hat{r}_{m}^{j}\mathbb{E}\Big(\frac{p_{m}^{j+1}
-2p_{m}^{j}+p_{m}^{j-1}}{\D t^{2}}\Big|\mathcal{F}_{j\D t}\Big)\Big]dx
=K\mathbb{E}\sum_{j=0}^{m-1}\int_{G}\frac{(\hat{r}_{m}^{j+1}-\hat{r}_{m}^{j})^{2}}{\D t^{2}}dx.
\end{equation}
Using Lemma \ref{lem-a1} again and combining with \eqref{2.11}, we conclude that
\begin{equation}\label{a14}
\begin{aligned}
&-\mathbb{E}\sum_{j=1}^{m-1}\int_{G}\mathbb{E}\Big(\frac{\hat{r}_{1m}^{j+1}-\hat{r}_{1m}^{j}}{\D t}
\Big|\mathcal{F}_{j\D t}\Big)\mathbb{E}\Big(\frac{p_{m}^{j+1}-2p_{m}^{j}+p_{m}^{j-1}}{\D t^{2}}
\Big|\mathcal{F}_{j\D t}\Big)dx \\
&=-\mathbb{E}\sum_{j=1}^{m-1}\int_{G}\mathbb{E}\Big(\frac{\hat{r}_{1m}^{j+1}
-\hat{r}_{1m}^{j}}{\D t}\Big|\mathcal{F}_{j\D t}\Big)\frac{1}{\D t}\mathbb{E}\Big(
\frac{p_{m}^{j+1}-p_{m}^{j}}{\D t}-\frac{p_{m}^{j}-p_{m}^{j-1}}{\D t}\Big|
\mathcal{F}_{j\D t}\Big)dx.
\end{aligned}
\end{equation}
Further simplification leads to
\begin{equation}\label{a15}
\begin{aligned}
& -\mathbb{E}\sum_{j=1}^{m-1}\int_{G}\mathbb{E}\Big(\frac{\hat{r}_{1m}^{j+1}
-\hat{r}_{1m}^{j}}{\D t} \Big| \mathcal{F}_{j\D t}\Big)\mathbb{E}\Big(\frac{p_{m}^{j+1}-2p_{m}^{j}+p_{m}^{j-1}}{\D t^{2}} \Big| \mathcal{F}_{j\D t}\Big)dx \\
&=\mathbb{E}\sum_{j=1}^{m-1}\int_{G}\frac{\rho^{K}}{\lambda^{2}}\mathbb{E}\Big(
\frac{\hat{r}_{1m}^{j+1}-\hat{r}_{1m}^{j}}{\D t} \Big| \mathcal{F}_{j\D t}\Big)\mathbb{E}\Big(\frac{\hat{r}_{1m}^{j+1}e^{-2\lambda\phi_{m}^{j+1}}
	-\hat{r}_{1m}^{j}e^{-2\lambda\phi_{m}^{j}}}{\D t} \Big| \mathcal{F}_{j\D t}\Big)dx\\
&=\mathbb{E}\sum_{j=1}^{m-1}\int_{G}\frac{\rho^{K}}{\lambda^{2}}\bigg\{\frac{[
	\mathbb{E}(\hat{r}_{1m}^{j+1}-\hat{r}_{1m}^{j}|\mathcal{F}_{j\D t})]^{2}}{\D t^{2}}
e^{-2\lambda\phi_{m}^{j}} \\
& \qq\qq\qq\q+\mathbb{E}\Big(\frac{\hat{r}_{1m}^{j+1}-\hat{r}_{1m}^{j}}{\D t} \Big| \mathcal{F}_{j\D t}\Big)\frac{e^{-2\lambda\phi_{m}^{j+1}}-e^{-2\lambda\phi_{m}^{j}}}{\D t}
\mathbb{E}(\hat{r}_{1m}^{j+1}|\mathcal{F}_{j\D t})\bigg\}dx,
\end{aligned}
\end{equation}
and
\begin{eqnarray}\label{a16} 
&& -\mathbb{E}\sum_{j=1}^{m-1}\int_{G}\lambda v_{m}^{j}e^{2\lambda\phi_{m}^{j}}
\mathbb{E}\Big(\frac{p_{m}^{j+1}-2p_{m}^{j}+p_{m}^{j-1}}{\D t^{2}}\Big|\mathcal{F}_{j\D t}\Big)dx \nonumber\\
& =&-\mathbb{E}\sum_{j=1}^{m-1}\int_{G}\lambda v_{m}^{j}e^{2\lambda\phi_{m}^{j}}\frac{1}{\D t}\mathbb{E}\Big(\frac{p_{m}^{j+1}-p_{m}^{j}}{\D t}-\frac{p_{m}^{j}-p_{m}^{j-1}}{\D t}
\Big|\mathcal{F}_{j\D t}\Big)dx\\
&=&\lambda\mathbb{E}\sum_{j=1}^{m-1}\int_{G}\frac{\rho^{K}}{\lambda^{2}}v_{m}^{j}
\Big[\mathbb{E}\Big(\frac{\hat{r}_{1m}^{j+1}-\hat{r}_{1m}^{j}}{\D t}e^{-2\lambda\phi_{m}^{j}}\Big|\mathcal{F}_{j\D t}\Big)
+\frac{e^{-2\lambda\phi_{m}^{j+1}}-e^{-2\lambda\phi_{m}^{j}}}{\D t}\mathbb{E}(\hat{r}_{1m}^{j+1}|\mathcal{F}_{j\D t})\Big]dx.\nonumber 
\end{eqnarray}
Combining equations (\ref{a11})-(\ref{a16}), we arrive at 
\begin{equation}\label{a17}
\begin{aligned}
&\mathbb{E}\sum_{j=0}^{m-1}\int_{G}\bigg\{\frac{(\hat{w}_{m}^{j+1}-\hat{w}_{m}^{j})^{2}}{\D t^{2}}e^{-2\lambda\phi_{m}^{j}}
+\frac{\rho^{K}}{\lambda^{2}}\frac{[\mathbb{E}(\hat{r}_{1m}^{j+1}-\hat{r}_{1m}^{j}|\mathcal{F}_{j\D t})]^{2}}{\D t^{2}}
e^{-2\lambda\phi_{m}^{j}}+K\frac{(\hat{r}_{m}^{j+1}-\hat{r}_{m}^j)^2}{\D t^2}\bigg\}dx\\
&=-\mathbb{E}\sum_{j=0}^{m-1}\int_{G}\frac{(\hat{w}_{m}^{j+1}-\hat{w}_{m}^{j})}{\D t}
\frac{(e^{-2\lambda\phi_{m}^{j+1}}-e^{-2\lambda\phi_{m}^{j}})}{\D t}\hat{w}_{m}^{j+1}dx\\
&\q-\mathbb{E}\sum_{j=1}^{m-1}\int_{G}\frac{\rho^{K}}{\lambda^{2}}\mathbb{E}
\Big(\frac{\hat{r}_{1m}^{j+1}-\hat{r}_{1m}^{j}}{\D t}\Big|\mathcal{F}_{j\D t}\Big)
\frac{e^{-2\lambda\phi_{m}^{j+1}}-e^{-2\lambda\phi_{m}^{j}}}{\D t}
\mathbb{E}(\hat{r}_{1m}^{j+1}|\mathcal{F}_{j\D t})dx \\
&\q-\lambda\mathbb{E}\sum_{j=1}^{m-1}\!\int_{G}\!\frac{\rho^{K}}{\lambda^{2}}v_{m}^{j}
\bigg[\mathbb{E}\Big(\frac{\hat{r}_{1m}^{j+1}\!-\!\hat{r}_{1m}^{j}}{\D t}e^{-2\lambda\phi_{m}^{j}}
\Big|\mathcal{F}_{j\D t}\Big)\!+\!\frac{(e^{-2\lambda\phi_{m}^{j+1}}\!-\!e^{-2\lambda\phi_{m}^{j}})}{\D t}
\mathbb{E}(\hat{r}_{1m}^{j+1}|\mathcal{F}_{j\D t})\bigg]dx.
\end{aligned}
\end{equation}
Applying H\"{o}lder's inequality,  there is a constant $C=C(K,\l)>0$, independent of $m$, such that
\begin{equation}\label{a18}
\begin{aligned}
&\mathbb{E}\sum_{j=0}^{m-1}\int_{G}\bigg\{\frac{(\hat{w}_{m}^{j+1}-\hat{w}_{m}^{j})^{2}}{\D t^{2}}
e^{-2\lambda\phi_{m}^{j}}+\frac{\rho^{K}}{\lambda^{2}}\frac{[\mathbb{E}
(\hat{r}_{1m}^{j+1}-\hat{r}_{1m}^{j}|\mathcal{F}_{j\D t})]^{2}}{\D t^{2}}e^{-2\lambda\phi_{m}^{j}} +K\frac{(\hat{r}_m^{j+1}-\hat{r}_m^j)^2}{\D t^2}\bigg\}dx \\
&\leq C\bigg[\mathbb{E}\sum_{j=1}^{m-1}\int_{G}(|\hat{w}_{m}^{j}|^{2}+|\hat{r}_{1m}^{j}|^{2}
+K|\hat{r}_{m}^{j}|^{2}+|v_{m}^{j}|^{2})dx+\dbE\int_{G}|\hat{r}_{1m}^{m}|^{2}dx\bigg].
\end{aligned}
\end{equation}
Finally, combining equations \eqref{2.13} and \eqref{a18}, we establish the desired estimate  \eqref{2.14}.

\ss

{\bf Step 5.} In this step, we prove \eqref{2.15}. We begin by substituting 
$\{w_{0m}^{j}\}_{j=0}^{m}$ with $\{p_{m}^{j}\}_{j=0}^{m}$ in equation \eqref{a5} and recalling that  $p_{m}^{j}=K\hat{r}_{m}^{j}$. This yields 
\begin{equation*}\label{a19}
\begin{aligned}
0&=-\mathbb{E}\!\sum_{j=1}^{m-1}\!\int_{G}\!\bigg\{p_{m}^{j}\Big[\mathbb{E}\Big(
\frac{p_{m}^{j+1}\!-\!2p_{m}^{j}\!+\!p_{m}^{j-1}}{\D t^{2}}\Big|\mathcal{F}_{j\D t}\Big)
-\sum_{j_{1},j_{2}=1}^{n}\!\partial_{x_{j_{2}}}(a^{j_{1}j_{2}}\partial_{x_{j_{1}}}p_{m}^{j})\Big]\! +\hat{w}_{m}^{j}p_{m}^{j}e^{-2\lambda\phi_{m}^{j}}\bigg\}dx\\
&=\mathbb{E}\!\sum_{j=1}^{m-1}\!\int_{G}\!\frac{(p_{m}^{j+1}-p_{m}^{j})^{2}}{\D t^{2}}dx+
\dbE\sum_{j=1}^{m-1}\!\int_{G}\!\sum_{j_{1},j_{2}=1}^{n}\!a^{j_{1}j_{2}}\partial_{x_{j_{1}}}p_{m}^{j}\partial_{x_{j_{2}}}p_{m}^{j}dx -\dbE\!\sum_{j=1}^{m-1}\!\int_{G}\!\hat{w}_{m}^{j}p_{m}^{j}e^{-2\lambda\phi_{m}^{j}}dx.
\end{aligned}
\end{equation*}
Combining this result with \eqref{2.13} and \eqref{2.14}, and recalling that $p_{m}^{j}=K\hat{r}_{m}^{j}$, we obtain the following uniform estimate: 
\begin{equation}\label{a20}
\D t\dbE\sum_{j=1}^{m-1} \int_{G}|\nabla\hat{r}_{m}^{j}|^{2}dx\leq C.
\end{equation}
Next, we multiply equation \eqref{2.5} by $\{\hat{w}_{m}^{j}\}_{j=0}^{m}$ to derive 
\begin{equation*}\label{a21}
\begin{aligned}
0&=-\mathbb{E}\sum_{j=1}^{m-1}\int_{G}\Big[\mathbb{E}\Big(\frac{\hat{w}_{m}^{j+1}-2\hat{w}_{m}^{j}+\hat{w}_{m}^{j-1}}{\D t^{2}} \Big| \mathcal{F}_{j\D t}\Big)-\sum_{j_{1},j_{2}=1}^{n}\partial_{x_{j_{2}}}(a^{j_{1}j_{2}}\partial_{x_{j_{1}}}\hat{w}_{m}^{j})\\
& \qq\qq\qq-\mathbb{E}\Big(\frac{\hat{r}_{1m}^{j+1}-\hat{r}_{1m}^{j}}{\D t} \Big| \mathcal{F}_{j\D t}\Big)-\lambda v_{m}^{j}e^{2\lambda\phi_{m}^{j}}-\hat{r}_{m}^{j} \Big]\hat{w}_{m}^{j}dx \\
&=\mathbb{E}\sum_{j=1}^{m-1}\int_{G}\frac{(\hat{w}_{m}^{j+1}-\hat{w}_{m}^{j})^{2}}{\D t^{2}}dx+
\dbE\sum_{j=1}^{m-1}\int_{G}\sum_{j_{1},j_{2}=1}^{n}a^{j_{1}j_{2}}\partial_{x_{j_{1}}}\hat{w}_{m}^{j}\partial_{x_{j_{2}}}\hat{w}_{m}^{j}dx\\
&\q+\dbE\sum_{j=1}^{m-1}\int_{G}\Big[\mathbb{E}\Big(\frac{\hat{r}_{1m}^{j+1}-\hat{r}_{1m}^{j}}{\D t} \Big| \mathcal{F}_{j\D t}\Big)+\lambda v_{m}^{j}e^{2\lambda\phi_{m}^{j}}+\hat{r}_{m}^{j} \Big]\hat{w}_{m}^{j}dx,
\end{aligned}
\end{equation*}
which together with \eqref{2.13}, \eqref{2.14} and \eqref{a20}, implies \eqref{2.15}. \endpf

\subsection{A  Carleman estimate for the hyperbolic operator $\cA$}\label{ssec2.2}

In this subsection, we establish the following Carleman estimate for the hyperbolic operator $\cA$. 

\begin{lemma}\label{lem3.1}
There exists $\ch{w}\in L^{2}_{\dbF}(\Omega;H_{0}^{1}(0,T;L^{2}(G)))\cap L^{2}_{\dbF}(\Omega;L^{2}(0,T;H^{1}_{0}(G)))$ satisfying
\begin{equation}\label{3.3}
\begin{cases}
\cA\ch{w}=\lambda\th^{2}v & \text{in}\q Q,\\
\ch{w}=0 &\text{on}\q  \Sigma,\\
\ch{w}(0)=\ch{w}(T)=0&\text{in}\q G.
\end{cases}
\end{equation}
Moreover, there exists $\l_{1}=\l_1 (\vp, T, (a^{jk})_{n\t n}, G, G_0)>0$ such that for $\l\geq\l_{1}$, it holds that
\begin{equation}\label{3.4}
\dbE\int_{Q}\th^{-2}(|\nabla \ch{w}|^{2}+\ch{w}_{t}^{2}+\l^{2}\ch{w}^{2})dxdt
\leq C\l\dbE\int_{Q}\th^{2}v^{2}dxdt.
\end{equation}
\end{lemma}

We begin by recalling the following Carleman estimate for the operator $\cA$.

\begin{lemma}\label{lem2.1}{\rm \cite[Theorem 2.1]{FLL23}}
Assume that Condition \ref{con1.1} holds and let $T>T_{0}$. Then there exist  positive constants $\mu_{0}=\mu_0(\vp, T, (a^{jk})_{n\t n}, G, G_0)$ and $C=C(\vp, T, (a^{jk})_{n\t n}, G, G_0 )$ such that for all $\mu\geq\mu_{0}$, there exists  $\l_{0}=\l_{0}(\mu, \vp, T, (a^{jk})_{n\t n}, G, G_0 )$, so that for all $\l\geq \l_{0}$, and for any $u\in H^{1}(Q)$ with $\cA u\in L^{2}(Q)$, $u|_{\Sigma}=0$ and $\supp u\subset[0,T]\t(\ol{G}\setminus G_{0})$,
it holds that
\begin{equation}\label{2.2}
\l\mu\int_{Q}\theta^{2}\( u_{t}^{2}+|\nabla u|^{2}+\lambda^{2}\mu^{2}
u^{2}\)dxdt\leq
C \int_{Q}\theta^{2}|\cA u|^{2}dxdt.
\end{equation}
\end{lemma}

\noindent{\it Proof of Lemma \ref{lem3.1}.}
We divide the proof into five steps.

{\bf Step 1.} Recalling the functions $\{(\hat{w}_{m}^{j},\hat{r}_{1m}^{j},\hat{r}_{m}^{j})\}_{j=0}^{m}$ from  Proposition \ref{prop2.1}, we define the following piecewise linear interpolations:
\begin{equation}\label{3.7}
\begin{cases} \tl{w}^{m}(t,x)\!\=\!\frac{1}{\D t}\sum_{j=0}^{m-1}\dbE\Big(\Big\{(t-j\D t)\hat{w}_{m}^{j+1}(x) -[t-(j+1)\D t]\hat{w}_{m}^{j}(x)\Big\}\Big|\cF_{j\D t}\Big)\chi_{(jh,(j+1)h]}(t),\\ \tl{r}^{m}_{1}(t,x)\!\=\!\hat{r}_{1m}^{0}\chi_{\{0\}}(t)\!+\!\frac{1}{\D t}\!\sum_{j=0}^{m-1}\!\dbE\Big(\!\Big\{(t\!-\!j\D t)\hat{r}_{1m}^{j+1}(x) \!-[t\!-(j\!+\!1)\D t]\hat{r}_{1m}^{j}(x)\Big\}\Big|\cF_{j\D t}\!\Big)\chi_{(j\D t,(j\!+\!1)\D t]}(t),\\ \tl{r}^{m}(t,x)\!\=\!\frac{1}{\D t}\sum_{j=0}^{m-1}\dbE\Big(\Big\{(t\!-\!j\D t)\hat{r}_{m}^{j+1}(x)
-[t-(j+1)\D t]\hat{r}_{m}^{j}(x)\Big\}\Big|\cF_{j\D t}\Big)\chi_{(j\D t,(j+1)\D t]}(t),
\end{cases}
\end{equation}
where $\chi$ denotes the characteristic function.
 
From \eqref{2.13}-\eqref{2.15}, we can extract a subsequence of  $\{(\tl{w}^{m},\tl{r}_{1}^{m},\tl{r}^{m})\}_{m=3}^{\i}$ that converges weakly to some $(\tl{w},\tl{r}_{1},\tl{r})\in (L_{\dbF}^{2}(\Omega;H^{1}(0,T;L^{2}(G))))^{3}$, as $m\rightarrow\i$. Moreover, by \eqref{3.5}, $\tl{w}\in L_{\dbF}^{2}(\Omega;H^{1}(0,T;L^{2}(G)))$ is the weak solution to the following random hyperbolic equation:
\begin{equation}\label{3.8}
\begin{cases}
\cA\tl{w}=\tl{r}_{1,t}+\lambda\th^{2}v+\tl{r} & \text{in}\q Q,\\
\tl{w}=0 &\text{on}\q {\Sigma},\\
\tl{w}(0)=\tl{w}(T)=0&\text{in}\q G.
\end{cases}
\end{equation}
Following the results in \cite[Appendix B]{FLLZ16}, we have the improved regularity:
\begin{equation}\label{3.9}
\tl{w}\in L_{\dbF}^{2}(\Omega;C([0,T];H^{1}_{0}(G)))\cap L_{\dbF}^{2}(\Omega;C^{1}([0,T];L^{2}(G))),
\end{equation}
For any $K>1$, define $\tl{p}\=K\tl{r}$. Proposition \ref{prop2.1} implies that $\tl{p}$ satisfies
\begin{equation}\label{3.10}
\begin{cases}
\cA\tl{p}+\th^{-2}\tl{w}=0 &\text{in}\q Q,\\
\tl{p}=0 &\text{on}\q {\Sigma},\\
\tl{p}(0)=\tl{p}(T)=0&\text{in}\q G,\\
\tl{p}_{t}+\rho^{K}\th^{-2}\frac{\tl{r}_{1}}{\l^{2}}=0&\text{in}\q Q.
\end{cases}
\end{equation}
Since $\tl{r}_{1}\in L_{\dbF}^{2}(\Omega;H^{1}(0,T;L^{2}(G)))$, we similarly obtain
\begin{equation}\label{3.11}
\tl{p}\in L_{\dbF}^{2}(\Omega;C([0,T];H^{1}_{0}(G)))\cap L_{\dbF}^{2}(\Omega;C^{1}([0,T];L^{2}(G))).
\end{equation}
By \eqref{3.10}, $\tl{p}_{t}$ satisfies
\begin{equation}\label{3.12}
\begin{cases}
\cA\tl{p}_{t}+(\th^{-2}\tl{w})_{t}=0 &\text{in}\q Q,\\
\tl{p}_{t}=0 &\text{on}\q \Sigma,\\
\tl{p}_{tt}+\frac{\rho^{K}}{\l}\th^{-2}\Big(\frac{\tl{r}_{1,t}}{\l}-2\phi_{t}
\tl{r}_{1}\Big)=0&\text{in}\q Q.
\end{cases}
\end{equation}
Applying Lemma \ref{lem2.1} to both $\tl{p}$ in \eqref{3.10} and $\tl{p}_{t}$ in \eqref{3.12}, and setting $\mu=\mu_0$, we conclude that there is a $\l_{1}>0$ such that for $\l\geq\l_{1}$, \vspace{-1mm}
\begin{equation}\label{3.13}
\lambda\dbE \int_{Q}\theta^{2}\(|\tl{p}_{t}|^{2}+|\nabla \tl{p}|^{2}+\lambda^{2}
|\tl{p}|^{2}\)dxdt\leq
C \dbE\int_{Q}\theta^{-2}\tl{w}^{2}dxdt,
\end{equation}
and\vspace{-1mm}
\begin{equation}\label{3.14}
\lambda\dbE\int_{Q}\theta^{2}\(|\tl{p}_{tt}|^{2}+|\nabla \tl{p}_{t}|^{2}+\lambda^{2}
|\tl{p}_{t}|^{2}\)dxdt
\leq
C\dbE\int_{Q}\theta^{-2}(\tl{w}_{t}^{2}+\l^{2}\tl{w})dxdt.
\end{equation}

{\bf Step 2.} By \eqref{3.10}, we have
\begin{equation*}\label{3.15}
-\dbE\int_{Q}\tl{r}_{1,t}\tl{p}dxdt
=\dbE\int_{Q}\tl{r}_{1}\tl{p}_{t}dxdt
=-\dbE\int_{Q}\th^{-2}\rho^{K}\frac{\tl{r}_{1}^{2}}{\l^{2}}dxdt.
\end{equation*}
This implies that
\begin{equation}\label{3.16}
\begin{aligned}
0
&=\dbE\lg\cA\tl{w}-\tl{r}_{1,t}-\l\th^{2}v-\tl{r},\tl{p}\rg_{L^2(Q)}\\
&=-\dbE\int_{Q}\th^{-2}\tl{w}^{2}dxdt-\dbE\int_{Q}\th^{-2}\rho^{K}\frac{\tl{r}_{1}^{2}}{\l^{2}}dxdt
-\l\dbE\int_{Q}\th^{2}v\tl{p}dxdt-K\dbE\int_{Q}\tl{r}^{2}dxdt.
\end{aligned}
\end{equation}
By \eqref{3.13}, \eqref{3.16} and Cauchy-Schwarz inequality, we obtain that
\begin{equation}\label{3.17}
\dbE\int_{Q}\th^{-2}\tl{w}^{2}dxdt+\dbE\int_{Q}\th^{-2}\rho^{K}\frac{\tl{r}_{1}^{2}}{\l^{2}}dxdt
+K\dbE\int_{Q}\tl{r}^{2}dxdt
\leq \frac{C}{\l}\dbE\int_{Q}\th^{2}v^{2}dxdt.
\end{equation}

{\bf Step 3.} From \eqref{3.10} and \eqref{3.12}, and observing that $\tl{p}_{tt}(0)=\tl{p}_{tt}(T)=0$ in $G$, we have
\begin{equation}\label{3.18}
\begin{aligned}
0
&=\dbE\lg\cA\tl{w}-\tl{r}_{1,t}-\l\th^{2}v-\tl{r},\tl{p}_{tt}\rg_{L^2(Q)}\\
&=-\dbE\int_{Q}\tl{w}(\th^{-2}\tl{w})_{tt}dxdt-\dbE\int_{Q}\tl{r}_{1,t}\tl{p}_{tt}dxdt
-\l\dbE\int_{Q}\th^{2}v\tl{p}_{tt}dxdt-\dbE\int_{Q}\tl{r}\tl{p}_{tt}dxdt.
\end{aligned}
\end{equation}
Recalling the definitions of 
$\th$ and $\phi$ from \eqref{3.2} and applying integration by parts, we obtain that\vspace{-1mm}
\begin{equation}\label{3.19}
-\dbE\int_{Q}\tl{w}(\th^{-2}\tl{w})_{tt}dxdt=
\dbE\int_{Q}\th^{-2}(\tl{w}_{t}^{2}+ \l\phi_{tt}\tl{w}^2-2\l^{2}\phi_{t}^2\tl{w}^2)dxdt.
\end{equation}
From \eqref{3.12}, we derive that\vspace{-1mm}
\begin{equation}\label{3.20}
-\dbE\int_{Q}\tl{r}_{1,t}\tl{p}_{tt}dxdt=
\dbE\int_{Q}\th^{-2}\rho^{K}\Big(\frac{\tl{r}^{2}_{1,t}}{\l^2}
-\frac{2}{\l}\phi_{t}\tl{r}_{1}\tl{r}_{1,t}\Big)dxdt.
\end{equation}
Using integration by parts and the relation $\tl{p}=K\tl{r}$, we get that
\begin{equation}\label{3.21}
-\dbE\int_{Q}\tl{r}\tl{p}_{tt}dxdt=K\dbE\int_{Q}\tl{r}_{t}^{2}dxdt.
\end{equation}
Combining equations \eqref{3.18}--\eqref{3.21} yields that
\begin{equation}\label{3.22}
\begin{aligned}
\l\dbE\int_{Q}\th^{2}v\tl{p}_{tt}dxdt
&=\dbE\int_{Q}\th^{-2}\rho^{K}\Big(\frac{\tl{r}^{2}_{1,t}}{\l^2}
-\frac{2}{\l}\phi_{t}\tl{r}_{1}\tl{r}_{1,t}\Big)dxdt\\
&\q+K\dbE\int_{Q}\tl{r}_{t}^{2}dxdt
+\dbE\int_{Q}\th^{-2}(\tl{w}_{t}^{2}+\l\phi_{tt}\tl{w}^2-2\l^{2}\phi_{t}^2\tl{w}^2)dxdt.
\end{aligned}
\end{equation}
Finally, applying \eqref{3.17} and the Cauchy-Schwarz inequality, along with \eqref{3.14}, we obtain that
\begin{equation}\label{3.23}
\dbE\int_{Q}\th^{-2}(\tl{w}_{t}^{2}+\l^{2}\tl{w}^{2})dxdt
+\dbE\int_{Q}\th^{-2}\rho^{K}\Big(\frac{\tl{r}_{1,t}^{2}}{\l^{2}}
+\tl{r}_{1}^{2}\Big)dxdt
\leq C\l\dbE\int_{Q}\th^{2}v^{2}dxdt.
\end{equation}

{\bf Step 4.} It follows from \eqref{3.8} that
\begin{eqnarray}\label{3.24}
&&\dbE\lg\tl{r}_{1,t}+\l \th^{2}v+\tl{r},\th^{-2}\tl{w}\rg_{L^2(Q)} \nonumber\\
&=&\dbE\lg\cA\tl{w},\th^{-2}\tl{w}\rg_{L^2(Q)}\nonumber\\
&=&-\dbE\int_{Q}\tl{w}_{t}(\th^{-2}\tl{w})_{t}dxdt+
\sum_{j,k=1}^{n}\dbE\int_{Q}a^{jk}(\th^{-2}\tl{w})_{x_{j}}\tl{w}_{x_{k}}dxdt\\
&=&-\dbE\int_{Q}\th^{-2}(\tl{w}_{t}^{2}+\l\phi_{tt}\tl{w}^2-2\l^{2}\phi_{t}^{2}\tl{w}^{2})dxdt
+\sum_{j,k=1}^{n}\dbE\int_{Q}\th^{-2}a^{jk}\tl{w}_{x_{j}}\tl{w}_{x_{k}}dxdt\nonumber\\
&&\q-2\l\sum_{j,k=1}^{n}\dbE\int_{Q}\th^{-2}a^{jk}\tl{w}_{x_{j}}\phi_{x_{k}} \tl{w}dxdt.\nonumber
\end{eqnarray}
This yields that\vspace{-1mm}
\begin{equation}\label{3.25}
\begin{aligned}
\dbE\int_{Q}\th^{-2}|\nabla \tl{w}|^{2}dxdt
&\leq C\dbE\int_{Q}\Big[\th^{-2}|\tl{r}_{1,t}+\tl{r}||\tl{w}|+\l|v\tl{w}|
+\th^{-2}\big(\tl{w}_{t}^{2}+\l^{2}\tl{w}^{2}\big)\Big]dxdt\\
&\leq C\dbE\int_{Q}\Big[\th^{2}v^{2}+\th^{-2}\Big(\frac{\tl{r}_{1,t}^2}{\l^2}+\tl{r}^{2}
+\tl{w}_{t}^{2}+\l^{2}\tl{w}^{2}\Big)\Big]dxdt.
\end{aligned}
\end{equation}
By combining estimates \eqref{3.17}, \eqref{3.23}, and \eqref{3.25}, and choosing the parameter 
$K$ to satisfy \vspace{-1mm}
\begin{equation}\label{3.26}
K\geq e^{2\l\max_{(t,x)\in{Q}}|\phi(t,x)|},
\end{equation}
we can eliminate the term $\dbE\int_{Q}\th^{-2}\tl{r}^{2}dxdt$ in \eqref{3.25}. Consequently, we obtain that\vspace{-1mm}
\begin{equation}\label{3.27}
\dbE\int_{Q}\th^{-2}(|\nabla \tl{w}|^{2}+\tl{w}_{t}^{2}+\l^{2}\tl{w}^{2})dxdt
+\dbE\int_{Q}\th^{-2}\rho^{K}\Big(\frac{\tl{r}_{1,t}^{2}}{\l^{2}}+\tl{r}_{1}^{2}\Big)dxdt
\leq C\l\dbE\int_{Q}\th^{2}v^{2}dxdt.
\end{equation}

{\bf Step 5.} Recall that the solution triple  $(\tl{w},\tl{r}_{1},\tl{r})$  depends on the parameter  $K$. To emphasize this dependence, we denote it by $(\tl{w}^{K},\tl{r}^{K}_{1},\tl{r}^{K})$.

Fixing $\l$ and letting $K\rightarrow \i$, we observe that $\rho^{K}(x)\rightarrow\i$  for all $x\notin G_{0}$. From  \eqref{3.17} and \eqref{3.27}, we conclude that there exists a subsequence of $\{(\tl{w}^{K},\tl{r}^{K}_{1},\tl{r}^{K})\}_{K=1}^{\i}$ that converges weakly to a limit $(\ch{w},0,0)$ in $
\big(L^{2}_{\dbF}(\Omega;H_{0}^{1}((0,T);L^{2}(G)))\cap L^{2}_{\dbF}(\Omega;L^{2}((0,T);H^{1}_{0}(G)))\big) 
\t L^{2}_{\dbF}(\Omega;H^{1}((0,T);L^{2}(G)))\t L^{2}_{\dbF}(\Omega;L^{2}(Q))$. 
Consequently, by equation \eqref{3.8}, the limit function $\ch{w}$ satisfies 
\begin{equation*}\label{3.28}
\begin{cases}
\cA\ch{w}=\lambda\th^{2}v & \text{in}\q Q,\\
\ch{w}=0 &\text{on}\q \Sigma,\\
\ch{w}(0)=\ch{w}(T)=0&\text{in}\q G.
\end{cases}
\end{equation*}
Applying estimate \eqref{3.25} to $\ch{w}$ yields the  estimate
\begin{equation*}\label{3.29}
\dbE\int_{Q}\th^{-2}(|\nabla \ch{w}|^{2}+\ch{w}_{t}^{2}+\l^{2}\ch{w}^{2})dxdt
\leq C\l\dbE\int_{Q}\th^{2}v^{2}dxdt,
\end{equation*}
which completes the proof. \endpf

\subsection{An $L^2$-Carleman estimate for the hyperbolic operator $\cA$}\label{ssec2.3}

In this subsection, we establish an $L^2$-Carleman estimate for the hyperbolic partial differential operator $\cA$.

\begin{theorem}\label{thm3.1}
Assume Condition \ref{con1.1} holds and set $\mu=\mu_{0}$ as in Lemma \ref{lem2.1}. Let $\alpha\in L^2_{\dbF}(\O;W^{1,\i}(0,T;$ $L^{\i}(G)))$ and suppose  $T>T_{0}$.
For any stochastic process $Y\in L_{\dbF}^{2}(\Omega;C([0,T];L^{2}(G)))\cap L_{\dbF}^{2}(0,T;H^{1}_0(G))$ satisfying the following conditions:
\begin{enumerate}
\item[(i)] $Y(0)=Y(T)=0$ in $G$;
\item[(ii)] $\supp Y\subset [0,T]\t(\ol{G}\setminus G_{0})$;
\item[(iii)] For any $\eta\in L_{\dbF}^{2}(\Omega;H_{0}^{1}(0,T;L^2(G)))\cap L_{\dbF}^{2}(\Omega;L^2(0,T;H_0^1(G)))$ with
$\cA\eta\in L_{\dbF}^{2}(\Omega;L^{2}(0,T;L^2(G)))$, it holds that
\begin{equation}\label{3.1}
\dbE\lg Y,\cA\eta-(\alpha\eta)_{t}\rg_{L^{2}(Q)}
=\dbE\lg F,\eta\rg_{H^{-1}(Q),H_{0}^{1}(Q)};
\end{equation}
\end{enumerate}
there exist positive constants $\l_2= \lambda_{2}( \vp, T, (a^{jk})_{n\t n}, G, G_0, \a)>0$ and $C=C( \vp, T, (a^{jk})_{n\t n}, G, G_0, \a)>0$ such that   for any $\lambda>\lambda_{2}$, the following Carleman estimate holds:
\begin{equation}\label{3.2}
\lambda\dbE\int_{Q}\th^{2}Y^{2}dxdt
\leq  C \dbE\|\th F\|^{2}_{H^{-1}(Q)}.
\end{equation}
\end{theorem}

\noindent{\it Proof.}
Substituting the test function  $\eta$ in \eqref{3.1} with $\ch{w}$ defined in \eqref{3.3}, we finds that
\begin{equation}\label{3.5}
\dbE\lg Y,\lambda\th^{2}Y-(\alpha\ch{w})_{t}\rg_{L^{2}(Q)}
=\dbE\lg F,\ch{w}\rg_{H^{-1}(Q),H_{0}^{1}(Q)}.
\end{equation}
This, together with Cauchy-Schwarz inequality, implies that
\begin{equation}\label{3.6}
\begin{aligned}
&\l\dbE\int_{Q}\th^{2}Y^{2}dxdt\\
&=\dbE\lg Y,(\alpha\ch{w})_{t}\rg_{L^{2}(Q)}+\dbE\lg F,\ch{w}\rg_{H^{-1}(Q),H_{0}^{1}(Q)}\\
&\leq
\frac{1}{2C_{1}}\Big(\dbE\|\th Y\|^2_{L^{2}(Q)}+\dbE\|\th F\|^{2}_{H^{-1}(Q)}\Big)
+C_{1}\Big(\dbE\|\th^{-1}(\alpha\ch{w})_{t}\|^2_{L^{2}(Q)}+ \dbE\|\th^{-1}\ch{w}\|^{2}_{H^{1}_{0}(Q)}\Big).
\end{aligned}
\end{equation}
Taking into account the bound \eqref{3.4} and selecting the constant $C_{1}<\frac{1}{4C}$, where $C$  is the constant appearing in \eqref{3.4}, we obtain the desired estimate \eqref{3.2}.\endpf

\begin{remark}
In recent years, Carleman estimates for stochastic hyperbolic equations have been extensively studied as a powerful tool for addressing inverse problems and unique continuation problems in stochastic hyperbolic systems(e.g.,\cite{Dou2025,L13,Wu2022,Yuan15,Z08}). While our primary application in this work focuses on establishing observability estimates for backward stochastic hyperbolic equations, we emphasize that Theorem \ref{thm3.1} possesses independent theoretical significance. Notably, these results may find important applications in studying inverse problems and unique continuation properties for backward stochastic hyperbolic equations. However, such investigations lie beyond the scope of the current paper and warrant separate consideration in future research.
\end{remark}

\section{Observability estimate for the equation \eqref{1.4}}\label{sec3}
 
This section presents the proof of Theorem \ref{thm1.2}, organized into two main parts. First, in Subsection \ref{ssec3.1}, we derive three crucial energy estimates for solutions to the backward stochastic hyperbolic equation \eqref{1.4}. These estimates will serve as fundamental tools for our subsequent analysis. Then, in Subsection \ref{ssec3.2}, we complete the proof of Theorem \ref{thm1.2} by  these preliminary results.

\subsection{Some energy estimates for the solutions to backward stochastic hyperbolic equations}\label{ssec3.1}

Let us define the energy functional
\begin{equation}\label{4.1}
\cE(t)\=\frac{1}{2}\big(\dbE\|z(t,\cdot)\|^{2}_{L^{2}(G)}+\dbE\|\hat{z}(t,\cdot)\|^{2}_{H^{-1}(G)}
\big),
\end{equation}
where the pair $(z,\hat{z})$ represents the first two components of the solution to the backward stochastic hyperbolic equation \eqref{1.4}.

We first prove the following energy estimate.

\begin{proposition}\label{prop4.1}
Let $0\leq s_{1}<s_{2}<t_{2}<t_{1}\leq T$. Then there exists a constant $C>0$ such that
\begin{equation}\label{4.2}
\int_{s_{2}}^{t_{2}}\cE(t)dt\leq C\Big(\dbE\int_{s_{1}}^{t_{1}}\|z\|^{2}_{L^{2}(G)}dt+\|a_{4}z+Z\|^{2}_{L^{2}_{\dbF}(0,T;L^{2}(G))}
+\|a_{5}\hat{z}+\h{Z}\|^{2}_{L^{2}_{\dbF}(0,T;H^{-1}(G))}\Big).
\end{equation}
\end{proposition}

\noindent{\it Proof.} We begin by selecting a smooth cut-off function $\psi\in C^{\infty}_{0}(0,T)$ with the following properties:
\begin{equation*}
\psi(t)=
\begin{cases}
1, & t\in[s_{2},t_{2}],\\
0, & t\in\big(0,\frac{s_{1}+s_{2}}{2}\big]\cup\big[\frac{t_{1}+t_{2}}{2},T\big).
\end{cases}
\end{equation*}

By the standard result of the existence of weak solutions to elliptic equations, the second-order differential operator
$$
\cL\=-\sum_{j,k=1}^{n}\partial_{x_{k}}(a^{jk}\partial_{x_{j}}) 
:H^{1}_{0}(G)\rightarrow H^{-1}(G)
$$
is an isomorphism between these Sobolev spaces.

Define $\zeta=\cL^{-1}z$. Applying It\^{o}'s formula to the system \eqref{1.4} yields 
\begin{eqnarray}\label{4.4} 
0&\3n=&\3n\dbE\int_{Q}\psi z\cL^{-1}\Big[d\hat{z}-\sum_{j,k=1}^{n}(a^{jk}z_{x_{j}})_{x_{k}}dt-(a_{1}z+a_{2}Z-a_{3}\h{Z})dt\Big]dx \nonumber\\
&\3n=&\3n\dbE\int_{Q}\Big\{-\psi_{t}z\cL^{-1}(dz)-\psi\hat{z}\cL^{-1}\hat{z}dt-\psi Z\cL^{-1}\h{Z}dt+\psi z^{2}dt-\psi z\cL^{-1}[(a_{1}-a_{2}a_{4})z]dt \\
&&\qq\qq-\psi z\cL^{-1}[a_{2}(a_{4}z+Z)]dt +\psi z\cL^{-1}[a_{3}(a_{5}\hat{z}+\h{Z})]dt-\psi z\cL^{-1}(a_{3}a_{5}\hat{z})dt\Big\}dx.\nonumber 
\end{eqnarray}
For the third term in \eqref{4.4}, we have the decomposition
\begin{equation}\label{4.5}
\begin{aligned}
-\dbE\int_{Q}\psi Z\cL^{-1}\h{Z}dtdx
&=-\dbE\int_{Q}\psi\Big[(a_{4}z+Z)\cL^{-1}(a_{5}\hat{z}+\h{Z})-(a_{4}z)\cL^{-1}(a_{5}\hat{z}+\h{Z})\\
&\qq\qq\q-(a_{4}z+Z)\cL^{-1}(a_{5}\hat{z})+(a_{4}z)\cL^{-1}(a_{5}\hat{z})\Big]dt.
\end{aligned}
\end{equation}
Using the relation $\zeta=\cL^{-1}z$, the first term in \eqref{4.4} becomes
\begin{eqnarray}\label{4.6} 
&&-\dbE\int_{Q}\psi_{t}z\cL^{-1}(dz)dx\nonumber\\
&=&\dbE\int_{Q} \sum_{j,k=1}^{n}\Big[(\psi_{t}a^{jk}\zeta_{x_{j}}d\zeta)_{x_{k}}+\frac{1}{2}d(\psi_{t}a^{jk}\zeta_{x_{j}}\zeta_{x_{j}})
-\frac{1}{2}\psi_{tt}a^{jk}\zeta_{x_{j}}\zeta_{x_{k}}dt-\frac{1}{2}\psi_{t}a^{jk}d\zeta_{x_{j}}d\zeta_{x_{k}}\Big] dx\nonumber\\
&=&-\frac{1}{2}\dbE\int_{Q}\Big[
\psi_{tt}\Big( \sum_{j,k=1}^{n}a^{jk}\zeta_{x_{j}}\zeta_{x_{k}}\Big)dt
+\psi_{t}\Big( \sum_{j,k=1}^{n}a^{jk}d\zeta_{x_{j}}d\zeta_{x_{k}}\Big)\Big]dx\\
&\leq& C\dbE\int_{0}^{T}\Big(\psi_{tt}\|z\|^{2}_{H^{-1}(G)}+\psi_{t}\|Z\|^{2}_{H^{-1}(G)}\Big)dt\nonumber\\
&\leq& C\dbE\int_{0}^{T}\Big[(\psi_{t}+\psi_{tt})\|z\|^{2}_{H^{-1}(G)}+\psi_{t}\|a_{4}z+Z\|^{2}_{H^{-1}(G)}\Big]dt.\nonumber
\end{eqnarray}

Combining estimates \eqref{4.4}-\eqref{4.6} and applying the Cauchy-Schwarz inequality along with the continuous embedding  $L^{2}(G)\hookrightarrow H^{-1}(G)$, we finally arrive at
\begin{equation*} 
\begin{aligned}
\dbE\int_{s_{2}}^{t_{2}}\|\hat{z}\|^{2}_{H^{-1}(G)}dt&\leq C\dbE\int_{Q}\psi\hat{z}\cL^{-1}\hat{z}dtdx\\
&\leq C\Big(\dbE\int_{s_{1}}^{t_{1}}\|z\|^{2}_{L^{2}(G)}dt+\|a_{4}z+Z\|^{2}_{L^{2}_{\dbF}(0,T;L^{2}(G))}
+\|a_{5}\hat{z}+\h{Z}\|^{2}_{L^{2}_{\dbF}(0,T;H^{-1}(G))}\Big).
\end{aligned}
\end{equation*}
\endpf

The second energy estimate is as follows:
\begin{proposition}\label{prop4.2}
There exists a positive constant $C>0$ such that for all $0\leq t\leq s\leq T$, the following estimate holds:
\begin{equation}\label{4.3}
\cE(t)\leq C\big(\cE(s)+\|a_{4}z+Z\|^{2}_{L^{2}_{\dbF}(0,T;L^{2}(G))}
+\|a_{5}\hat{z}+\h{Z}\|^{2}_{L^{2}_{\dbF}(0,T;H^{-1}(G))}\big).
\end{equation}
\end{proposition}

To prove Proposition \ref{prop4.2}, we consider the following  random hyperbolic equation:
\begin{equation}\label{4.8}
\begin{cases}
\beta_{tt}-\sum_{j,k=1}^{n}(a^{jk}\beta_{x_{j}})_{x_{k}} =a_{1}\beta &\text{in}\quad (\h{T},T)\times G,\\
\beta=0 &\text{on}\quad (\h{T},T)\times \Gamma,\\
\beta(\h{T})=\beta_{0},\hat{\beta}(\h{T})=\beta_{1} &\text{in}\quad  G,
\end{cases}
\end{equation}
where $\h{T}\in[0,T)$.  For any initial data $(\beta_{0},\beta_{1})\in L^{2}_{\cF_{\h{T}}}(\Omega;H^{1}_{0}(G)\times L^{2}(G))$, the system \eqref{4.8} admits a unique solution
$$
\beta\in L_{\dbF}^{2}(\Omega;C([\h{T},T];H^{1}_{0}(G)))\cap L_{\dbF}^{2}(\Omega;C^{1}([\h{T},T];L^{2}(G))).
$$

Define the energy functional for this system by
\begin{equation}\label{4.9}
E_{1}(t)=\frac{1}{2}\dbE\int_{G}(|\beta(t,x)|^{2}+|\nabla \beta(t,x)|^{2}+|\beta_{t}(t,x)|^{2})dx,
\quad t\in [\h{T},T].
\end{equation}
Then we have the following result.
\begin{lemma}\label{lem4.1}
There exists a constant $C>0$ such that for any solution to \eqref{4.8} and for all  $\h{T}\leq t\leq s\leq T$, we have
\begin{equation}\label{4.10}
E_{1}(s)\leq C E_{1}(t).
\end{equation}
\end{lemma}

\noindent{\it Proof.} Let us define the modified energy functional
\begin{equation}\label{4.11}
\wt{E}_{1}(t)=\frac{1}{2}\dbE\int_{G}\Big(|\beta(t,x)|^{2}+\sum_{j,k=1}^{n}a^{jk}\beta_{x_{j}}(t,x) \beta_{x_{k}}(t,x)+|\beta_{t}(t,x)|^{2}\Big)dx,
\quad t\in [\h{T},T].
\end{equation}
Applying It\^o's formula to this energy functional yields 
\begin{equation*}\label{4.12}
\begin{aligned}
\wt{E}_{1}(s)-\wt{E}_{1}(t)&=\dbE\int_{t}^{s}\int_{G}\Big(\beta\beta_{t}+\sum_{j,k=1}^{n}a^{jk}\beta_{x_{j}}\beta_{tx_{k}}
+\beta_{t}\sum_{j,k=1}^{n}(a^{jk}\beta_{x_{j}})_{x_{k}}+a_{1}\beta\beta_{t}\Big) dxd\tau\\
&=\dbE\int_{t}^{s}\int_{G}(1+a_{1})\beta\beta_{t}dxd\tau,
\end{aligned}
\end{equation*}
where the cancellation of terms follows from integration by parts and the boundary conditions. This leads to the inequality
\begin{equation}\label{4.13}
\wt{E}_{1}(s)\leq \wt{E}_{1}(t)+C\int_{t}^{s}\wt{E}_{1}(\tau)d\tau.
\end{equation}
An immediate application of Gronwall's inequality to \eqref{4.13} establishes the desired energy estimate \eqref{4.10}. \endpf
\ms

\noindent{\it Proof of Proposition \ref{prop4.2}.}
Let $\dbS$ denote the unit sphere in the space $L^{2}_{\cF_{t}}(\Omega;H^{1}_{0}(G))\times L^{2}_{\cF_{t}}(\Omega;L^{2}(G))$. Taking $\h{T}=t$ in \eqref{4.8} and applying It\^o's formula, we obtain that
\begin{equation*}\label{4.28}
\begin{aligned}
&\dbE\lg \hat{z}(s),\beta(s)\rg_{H^{-1}(G),H^{1}_{0}(G)}-\dbE\lg \hat{z}(t),\beta_{0}\rg_{H^{-1}(G),H^{1}_{0}(G)}\\
&=\dbE\lg z(s),\hat{\beta}(s)\rg_{L^{2}(G)}-\dbE\lg z(t),\beta_{1}\rg_{L^{2}(G)}\\
&\q+\dbE\int_{t}^{s}\int_{G}\Big[
a_{2}\beta(a_{4}z+Z)-a_{3}\beta(a_{5}\hat{z}+\h{Z})
-\beta(a_{2}a_{4}z-a_{3}a_{5}\hat{z})\Big]dxd\tau.
\end{aligned}
\end{equation*}
From this identity, we derive the following energy estimate
\begin{eqnarray*}  
2\cE(t)
&=&\sup_{(\beta_{0},\beta_{1})\in\dbS}\Big|\dbE\lg z(t),\beta_{1}\rg_{L^{2}(G)}+\dbE\lg \hat{z}(t),-\beta_{0}\rg_{H^{-1}(G),H^{1}_{0}(G)}\Big|^{2}\\
&=&\sup_{(\beta_{0},\beta_{1})\in\dbS}\Big|\dbE\lg z(s),\beta(s)\rg_{L^{2}(G)}+\dbE\lg \hat{z}(s),-\hat{\beta}(s)\rg_{H^{-1}(G),H^{1}_{0}(G)}\\
&&\qq\qq\q+\dbE\int_{s}^{t}\int_{G}\Big[a_{2}\beta(a_{4}z+Z)-a_{3}\beta(a_{5}\hat{z}+\h{Z})
-\beta(a_{2}a_{4}z-a_{3}a_{5}\hat{z})\Big]d\tau\Big|^{2}\\
&\leq& C\Big[\cE(s)\sup_{(\beta_{0},\beta_{1})\in\dbS}\|(\beta(s),\hat{\beta}(s))\|^{2}_{L^{2}_{\cF_{s}}
(\Omega;H^{1}_{0}(G)\times L^{2}(G))}\\
&&\qq\q+\sup_{(\beta_{0},\beta_{1})\in\dbS}\|\beta\|^{2}_{L^{2}_{\dbF}(s,t;H^{1}_{0}(G))}\big(\|a_{4}z+Z\|^{2}_{L^{2}_{\dbF}(s,t;L^{2}(G))}
+\|a_{5}\hat{z}+\h{Z}\|^{2}_{L^{2}_{\dbF}(s,t;H^{-1}(G))}\big)\\
&&\qq\q+\sup_{(\beta_{0},\beta_{1})\in\dbS}\|\beta\|^{2}_{L^{2}_{\dbF}(s,t;H^{1}_{0}(G))}\int_{s}^{t}\cE(\tau)d\tau
\Big].
\end{eqnarray*}
This, together with \eqref{4.10} and Gronwall's inequality, implies \eqref{4.3}. \endpf

Finally, we give the third energy estimate.
\begin{proposition}\label{prop4.3}
There exists a constant $C>0$ such that for all $0\leq s\leq t\leq T$, the following estimate holds:
\begin{equation}\label{4.3-1}
\cE(t)\leq C\big(\cE(s)+\|a_{4}z+Z\|^{2}_{L^{2}_{\dbF}(0,T;L^{2}(G))}
+\|a_{5}\hat{z}+\h{Z}\|^{2}_{L^{2}_{\dbF}(0,T;H^{-1}(G))}\big).
\end{equation}
\end{proposition}

To prove this result, we consider the backward stochastic hyperbolic system:
\begin{equation}\label{4.14}
\begin{cases}
d\By=\hat{\By}dt+\BY dW(t) &\text{in}\quad (0,\h{T}')\times G,\\
d\hat{\By}-\sum_{j,k=1}^{n}(a^{jk}\By_{x_{j}})_{x_{k}} dt=a_{1}\By dt+\h{\BY} dW(t) &\text{in}\quad (0,\h{T}')\times G,\\
\By=0 &\text{on}\quad (0,\h{T}')\times \Gamma,\\
\By(\h{T}')=\By_{0},\hat{\By}(\h{T}')=\By_{1} &\text{in}\quad  G,
\end{cases}
\end{equation}
where $\h{T}'\in(0,T]$. According to \cite[Theorem 4.10]{LZ21}, for any initial data  $(\By_{0},\By_{1})\in L^{2}_{\cF_{\h{T}'}}(\Omega;H^{1}_{0}(G)\times L^{2}(G))$, the system \eqref{4.14} admits a unique weak solution $(\By,\hat{\By},\BY,\h{\BY})$ in the space
$
L^{2}_{\dbF}(\Omega;C([0,\h{T}'];H^{1}_{0}(G)))$ $\times
L^{2}_{\dbF}(\Omega;C([0,\h{T}'];L^{2}(G)))
\times L^{2}_{\dbF}(0,\h{T}';H^{1}_{0}(G))\times
L^{2}_{\dbF}(0,\h{T}';L^{2}(G))$.

Define the energy functional:
\begin{equation}\label{4.17}
E_{2}(t)=\frac{1}{2}\dbE\int_{G}(\By(t,x)^{2}+|\nabla \By(t,x)|^{2}+\hat{\By}(t,x)^{2})dx,
\quad t\in [0,\h{T}'].
\end{equation}
Then we have the following energy estimates:
\begin{lemma}\label{lem4.3}
There exists a constant $C>0$ such that for any solution to \eqref{4.14} and for all $0\leq s\leq t\leq \h{T}'$, the following estimates hold:
\begin{equation}\label{4.18}
E_{2}(s)\leq C E_{2}(t),
\end{equation}
and
\begin{equation}\label{4.19}
\dbE\int_{0}^{\h{T}'}\int_{G}(|\BY|^{2}+|\nabla\BY|^{2}+|\h{\BY}|^{2})dxdt
\leq C\|(\By_{0},\By_{1})\|^{2}_{L^{2}_{\cF_{\h{T}'}}(\Omega;H^{1}_{0}(G)\t L^{2}(G))}.
\end{equation}
\end{lemma}

\noindent{\it Proof.} We introduce the modified energy functional:
\begin{equation*}\label{4.20}
\wt{E}_{2}(t)=\frac{1}{2}\dbE\int_{G}\Big(|\By(t,x)|^{2}+\sum_{j,k=1}^{n}a^{jk}\By_{x_{j}}(t,x) \By_{x_{k}}(t,x)+|\hat{\By}(t,x)|^{2}\Big)dx,
\quad t\in [0,\h{T}'].
\end{equation*}
Applying It\^o's formula yields 
\begin{eqnarray}\label{4.21} 
\wt{E}_{2}(t)-\wt{E}_{2}(s)&\3n=\3n&\dbE\int_{s}^{t}\int_{G}\Big[\By\hat{\By}+\sum_{j,k=1}^{n}a^{jk}\By_{x_{j}}\hat{\By}_{x_{k}}+
\hat{\By}\sum_{j,k=1}^{n}(a^{jk}\By_{x_{j}})_{x_{k}}+a_{1}\By\hat{\By}\nonumber\\
&&\qquad\qquad\qquad+\frac{1}{2}\Big(\BY^{2}+\sum_{j,k=1}^{n}a^{jk}\BY_{x_{j}}\BY_{x_{k}}
+\h{\BY}^{2}\Big)\Big] dxd\tau\\
&\3n=\3n&\dbE\int_{s}^{t}\int_{G}\Big[(1+a_{1})\By\hat{\By}+\frac{1}{2}\BY^{2}+\frac{1}{2}\sum_{j,k=1}^{n}a^{jk}\BY_{x_{j}}\BY_{x_{k}}
+\frac{1}{2}\h{\BY}^{2}\Big]dxd\tau.\nonumber
\end{eqnarray}
This yields that
\begin{equation}\label{4.22}
\frac{1}{2}\dbE\int_{s}^{t}\int_{G}\Big(\BY^{2}+\sum_{j,k=1}^{n}a^{jk}\BY_{x_{j}}\BY_{x_{k}}
+\h{\BY}^{2}\Big)dxd\tau
+\wt{E}_{2}(s)\leq \wt{E}_{2}(t)+C\int_{s}^{t}\wt{E}_{2}(\tau)d\tau.
\end{equation}
From \eqref{4.22} and Gronwall's inequality, we obtain that
\begin{equation}\label{4.23}
\wt{E}_{2}(s)\leq C\wt{E}_{2}(t),
\end{equation}
which implies \eqref{4.18}. 

Combining \eqref{4.21} and \eqref{4.23}, we obtain the inequality \eqref{4.19}.
\endpf
 
\ss

\noindent{\it Proof of Proposition \ref{prop4.3}.}
Let $\h{T}'=t$ in \eqref{4.14}. It follows from It\^{o}'s formula that
\begin{equation}\label{4.24}
\begin{aligned}
&\dbE\lg \hat{z}(t),\By_{0}\rg_{H^{-1}(G),H^{1}_{0}(G)}-\dbE\lg \hat{z}(s),\By(s)\rg_{H^{-1}(G),H^{1}_{0}(G)}\\
&=\dbE\lg z(t),\By_{1}\rg_{L^{2}(G)}-\dbE\lg z(s),\hat{\By}(s)\rg_{L^{2}(G)}\\
&\q+\dbE\int_{s}^{t}\int_{G}\Big[a_{2}\By(a_{4}z+Z)-a_{3}\By(a_{5}\hat{z}+\h{Z})
-\By(a_{2}a_{4}z-a_{3}a_{5}\hat{z})\\
&\qq\qq\qq  +\BY(a_{5}\hat{z}+\h{Z})-\h{\BY}(a_{4}z+Z) -a_{5}\hat{z}\BY+a_{4}z\h{\BY}\Big]d\tau.
\end{aligned}
\end{equation}
It is easy to see that
\begin{equation}\label{4.25}
\begin{aligned}
&\Big|\dbE\int_{s}^{t}\int_{G}\Big[a_{2}\By(a_{4}z+Z)-a_{3}\By(a_{5}\hat{z}+\h{Z})
-\By(a_{2}a_{4}z-a_{3}a_{5}\hat{z})\Big]d\tau\Big|\\
&\leq C\Big[ \|\By\|_{L^{2}_{\dbF}(s,t;H^{1}_{0}(G))}\big(\|a_{4}z+Z\|_{L^{2}_{\dbF}(s,t;L^{2}(G))}
+\|a_{5}\hat{z}+\h{Z}\|_{L^{2}_{\dbF}(s,t;H^{-1}(G))}\big)\\
&\qq+\|\By\|_{L^{2}_{\dbF}(s,t;H^{1}_{0}(G))}\Big(\int_{s}^{t}\cE(\tau)d\tau\Big)^{\frac{1}{2}}\Big],
\end{aligned}
\end{equation}
and
\begin{equation}\label{4.26}
\begin{aligned}
&\Big|\dbE\int_{s}^{t}\int_{G}\Big[\BY(a_{5}\hat{z}+\h{Z})-\h{\BY}(a_{4}z+Z) -a_{5}\hat{z}\BY+a_{4}z\h{\BY}\Big]d\tau\Big|\\
&\leq C\Big[\|(\BY,\h{\BY})\|_{L^{2}_{\dbF}(s,t;L^{2}(G)\t H^{-1}(G))}\big(\|a_{4}z+Z\|_{L^{2}_{\dbF}(s,t;L^{2}(G))}
+\|a_{5}\hat{z}+\h{Z}\|_{L^{2}_{\dbF}(s,t;H^{-1}(G))}\big)\\
&\qq+\|(\BY,\h{\BY})\|_{L^{2}_{\dbF}(s,t;L^{2}(G)\t H^{-1}(G))}\Big(\int_{s}^{t}\cE(\tau)d\tau\Big)^{\frac{1}{2}}\Big].
\end{aligned}
\end{equation}
The energy functional $\cE(t)$ can be expressed as
\begin{eqnarray*}\label{4.27} 
2\cE(t)
&\3n=\3n&\sup_{(\By_{0},\By_{1})\in\dbS}\Big|\dbE\lg z(t),\By_{1}\rg_{L^{2}(G)}+\dbE\lg \hat{z}(t),-\By_{0}\rg_{H^{-1}(G),H^{1}_{0}(G)}\Big|^{2}\\
&\3n=\3n&\sup_{(\By_{0},\By_{1})\in\dbS}\Big|\dbE\lg z(s),\By(s)\rg_{L^{2}(G)}+\dbE\lg \hat{z}(s),-\hat{\By}(s)\rg_{H^{-1}(G),H^{1}_{0}(G)}\\
&&\qq\qq\q+\dbE\int_{s}^{t}\int_{G}\Big[a_{2}\By(a_{4}z+Z)-a_{3}\By(a_{5}\hat{z}+\h{Z})
-\By(a_{2}a_{4}z-a_{3}a_{5}\hat{z})\Big]d\tau\\
&&\qq\qq\q+\dbE\int_{s}^{t}\int_{G}\Big[\BY(a_{5}\hat{z}+\h{Z})-\h{\BY}(a_{4}z+Z) -a_{5}\hat{z}\BY+a_{4}z\h{\BY}\Big]d\tau\Big|^{2}.
\end{eqnarray*}
Combining the estimates \eqref{4.18}, \eqref{4.19}, \eqref{4.24}-\eqref{4.26} with Gronwall's inequality, we obtain the desired inequality \eqref{4.3-1}. \endpf

\subsection{Proof of Theorem \ref{thm1.2}}\label{ssec3.2}

\noindent{\it Proof of Theorem \ref{thm1.2}.} The proof proceeds in three main steps. 

{\bf Step 1}. In this step we choose suitable cut-off functions. We begin by defining the time intervals and key parameters as follows:\vspace{-2mm}
$$
T_{j}=\frac{T}{2}-\vep_{j}T,\quad
T_{j}'=\frac{T}{2}+\vep_{j}T,\q
R_{0}=\min_{x\in \ol{G}}\sqrt{\vp(x)},\q
R_{1}=\max_{x\in \ol{G}}\sqrt{\vp(x)},\vspace{-2mm}
$$
where $j=1,2$ and $0<\vep_{1}<\vep_{2}<\frac{1}{2}$.  
Since $c_{1}>\frac{\sqrt{M_{1}}}{T}>\frac{M_{1}}{T^{2}}$, from \eqref{1.9-1} and \eqref{2.1},  we have that for any $T>T_{0}$,  \vspace{-2mm}
\begin{equation}\label{5.1}
\si(0,x)=\si(T,x)\leq R_{1}^{2}-\frac{c_{1}T^{2}}{4}
<R_{1}^{2}-\frac{M_{1}}{4}<0,\quad \forall x\in G.
\end{equation}
Let $J=(0,T_{2})\cup(T_{2}',T)$. There exists $\vep_{2}$ sufficiently close to  $\frac{1}{2}$ such that\vspace{-2mm}
\begin{equation}\label{5.2}
\si(t,x)\leq \frac{R_{1}^{2}}{2}-\frac{c_{1}T^{2}}{8}<0,\quad
\forall(t,x)\in J\times G.
\end{equation}
Moreover, since\vspace{-2mm}
$$
\si\Big(\frac{T}{2},x\Big)=\vp(x)\geq R_{0}^{2},\quad \forall x\in G,\vspace{-2mm}
$$
there exists $\vep_{1}$ close to $0$ such that\vspace{-2mm}
\begin{equation}\label{5.3}
\si(t,x)\geq \frac{R_{0}^{2}}{2},\quad \forall(t,x)\in(T_{1},T_{1}')\times G.
\end{equation}
We now select    a cut-off function $\kappa_{1}\in C^{\infty}_{0}(0,T)$ such that\vspace{-2mm}
\begin{equation}\label{5.4}
\kappa_{1}(t)=1 \quad\text{in} \quad (T_{2},T_{2}'),
\end{equation}
and a cut-off function $\k_{2}\in C_{0}^{\infty}(\dbR^{n})$ satisfying that\vspace{-2mm}
\begin{equation}\label{5.6}
\k_{2}(x)=
\begin{cases}
1,& x\in \ol{G}\setminus G_{0},\\
0,& x\in \cO_{\delta/2}(\Gamma_{0})\cap G.
\end{cases}
\end{equation}

{\bf Step 2.} In this step, we prove that there is a $\lambda_{3}>0$ such that for $\lambda\geq\lambda_{3}$, it holds that
\begin{equation}\label{5.7}
\begin{aligned}
\lambda\dbE\int_{Q}e^{2\lambda\phi}z^{2}dxdt
\leq & C\Big(\|e^{\lambda\phi}(a_{4}z+Z)\|_{L_{\dbF}^{2}(0,T;L^{2}(G))}^{2}
+\|e^{\lambda\phi}(a_{5}\hat{z}+\h{Z})\|^{2}_{L_{\dbF}^{2}(0,T;H^{-1}(G))}\\
&\quad+\dbE\|z\|^{2}_{L^{2}(J\times G)}+\lambda^{2}\dbE\int_{0}^{T}\int_{G_{0}}e^{2\lambda\phi}z^{2}dxdt\Big).
\end{aligned}
\end{equation}

Define $\xi=\k_{1}\k_{2} z$. Then $\xi$ satisfies
\begin{equation}\label{5.8}
\begin{cases}
d\xi=\hat{\xi}dt+\k_{1}\k_{2} ZdW(t) &\text{in}\quad Q,\\
d\hat{\xi}-\sum_{j,k=1}^{n}(a^{jk}\xi_{x_{j}})_{x_{k}} dt=Fdt+\k_{2}(\k_{1}\h{Z}+\k_{1,t}Z)dW(t) & \text{in}\quad Q,\\
\xi=0 &\text{on}\quad \Sigma,\\
\xi(0)=\xi(T)=0 & \text{in}\quad G,
\end{cases}
\end{equation}
where $\hat{\xi}=\k_{2}(\k_{1}\hat{z}+\k_{1,t}z)$, and
\begin{equation}\label{5.9}
\begin{aligned}
F&=(\k_{1,tt}+a_{1}\k_{1}-a_{2}a_{4}\k_{1})\k_{2} z+(2\k_{1,t}+a_{3}a_{5}\k_{1})\k_{2}\hat{z}+\k_{1}\k_{2}a_{2}(a_{4}z+Z)\\
&\q-\k_{1}\k_{2} a_{3}(a_{5}\hat{z}+\h{Z}) -\k_{1}\sum_{j,k=1}^{n}(a^{jk}\k_{2,x_{j}}z)_{x_{k}}
-\k_{1}\sum_{j,k=1}^{n}a^{jk}\k_{2,x_{k}}z_{x_{j}}.
\end{aligned}
\end{equation}

For any $\eta\in L_{\dbF}^{2}(\Omega;H_{0}^{1}(Q))\=L_{\dbF}^{2}(\Omega;H_{0}^{1}(0,T;L^2(G))\cap L_{\dbF}^{2}(0,T;H_{0}^{1}(G))$ with $\cA\eta\in L_{\dbF}^{2}(0,T;L^{2}(G))$,
we have
\begin{equation}
	\ba{ll}
\ds \eta \[d\hat \xi-\sum_{j,k=1}^n (a^{jk} \xi_{x_j})_{x_k}dt \]\\
\ns\ds =d(\hat \xi\eta )-d(\xi \eta_t )-\sum_{j,k=1}^n (a^{jk}\xi_{x_j}\eta)_{x_k}+\sum_{j,k=1}^n (a^{jk}\eta_{x_j}\xi)_{x_k} -\xi \cA \eta \\
\ns\ds = \eta Fdt+\eta\k_2 (\k_1 \hat Z+\k_{1,t}Z)dW(t).
\ea 
\end{equation}
Immediately, we get
$$
\dbE\lg \xi,\cA\eta+(a_{3}a_{5}\eta)_{t}\rg_{L^{2}(Q)}=\dbE\lg F',\eta\rg_{H^{-1}(Q),H_{0}^{1}(Q)},
$$
where
\begin{equation}\label{5.10}
\begin{aligned}
F'&=(\k_{1,tt}+a_{1}\k_{1}-a_{2}a_{4}\k_{1}-a_{3}a_{5}\k_{1,t})\k_{2} z+2\k_{1,t}\k_{2}\hat{z}\\
&\q+\k_{1}\k_{2}a_{2}(a_{4}z+Z)-\k_{1}\k_{2} a_{3}(a_{5}\hat{z}+\h{Z})-\k_{1}\sum_{j,k=1}^{n}(a^{jk}\k_{2,x_{j}}z)_{x_{k}}
-\k_{1}\sum_{j,k=1}^{n}a^{jk}\k_{2,x_{k}}z_{x_{j}}.
\end{aligned}
\end{equation}
Taking $\alpha=-a_{3}a_{5}$ in Theorem \ref{thm3.1},
then for $\lambda\geq\lambda_{2}$, it follows from \eqref{3.2} that
\begin{equation}\label{5.11}
\lambda\dbE\int_{Q}\th^{2}\xi^{2}dxdt
\leq C \dbE\|\th F'\|^{2}_{H^{-1}(Q)}.
\end{equation}
For   $\vartheta\in L_{\dbF}^{2}(\Omega;H_{0}^{1}(Q))$, using integration by parts, 
we get that
\begin{equation}\label{5.12}
\begin{aligned}
&\dbE\lg \th F',\vartheta\rg_{H^{-1}(Q),H^{1}_{0}(Q)}\\
&=\dbE\int_{0}^{T}\int_{G}\th\Big[\k_{1}\k_{2} a_{2}(a_{4}z+Z)\vartheta-\k_{1}\k_{2} a_{3}(a_{5}\hat{z}+\h{Z})\vartheta-a_{2}a_{4}\k_{1}\k_{2} \vartheta z\\
&\qq\qq\qq-(\k_{1,tt}\vartheta+2\l\phi_{t}\k_{1,t}\vartheta+2\k_{1,t}\vartheta_{t}
+a_{3}a_{5}\k_{1,t}\vartheta-a_{1}\k_{1}\vartheta)\k_{2}z\\
&\qq\qq\qq+\k_{1}\Big(2\sum_{j,k=1}^{n}a^{jk}\k_{2,x_{j}}\vartheta_{x_{k}}
+\sum_{j,k=1}^{n}(a^{jk}\k_{2,x_{k}})_{x_{j}}\vartheta
+2\l\sum_{j,k=1}^{n}a^{jk}\phi_{x_{j}}\k_{2,x_{k}}\vartheta
\Big)z\Big]dxdt.
\end{aligned}
\end{equation}
This yields that
\begin{eqnarray}\label{5.13}
\dbE\|\th F'\|^{2}_{H^{-1}(Q)}
&\3n=\3n&\sup_{\|\vartheta\|_{L^{2}_{\dbF}(\Omega;H^{1}_{0}(Q))}=1}\Big|\dbE\lg \th F',\vartheta\rg_{H^{-1}(Q),H^{1}_{0}(Q)}\Big|^{2}\nonumber\\
&\3n\leq\3n& C\Big[ \|e^{\lambda\phi}(a_{4}z+Z)\|_{L_{\dbF}^{2}(0,T;L^{2}(G))}^{2}+
\|e^{\lambda\phi}(a_{5}\hat{z}+\h{Z})\|^{2}_{L_{\dbF}^{2}(0,T;H^{-1}(G))}
+\dbE\|e^{\lambda\phi}z\|^{2}_{L^{2}(Q)}\nonumber\\
&&\qq+\lambda^{2}e^{2\lambda e^{\mu_0(R_1^2-c_1 T^{2}/4)/2}}\|z\|^{2}_{L^{2}_{\dbF}(J;L^{2}(G))}
+\l^{2}\dbE\int_{0}^{T}\int_{G_{0}}e^{2\l\phi}z^{2}dxdt\Big].
\end{eqnarray}
Furthermore, by \eqref{5.2}, \eqref{5.4} and \eqref{5.6}, we have
\begin{eqnarray}\label{5.14} 
&&\dbE\int_{Q}\th^{2}\xi^{2}dxdt \nonumber\\
&=&\|e^{\lambda\phi}z\|^{2}_{L^{2}_{\dbF}(0,T;L^{2}(G))}-\dbE\int_{Q}e^{2\lambda\phi}
(1-\k_{2}^{2})z^{2}dxdt-\dbE\int_{Q}e^{2\lambda\phi}
(1-\k_{1}^{2})\k_{2}^{2}z^{2}dxdt\\
&\geq& \|e^{\lambda\phi}z\|^{2}_{L^{2}_{\dbF}(0,T;L^{2}(G))}
-\dbE\int_{0}^{T}\int_{G_{0}}e^{2\lambda\phi}z^{2}dxdt
-Ce^{2\lambda e^{\mu_0(R_1^2-c_1 T^{2}/4)/2}}
\|z\|^{2}_{L^{2}_{\dbF}(J;L^{2}(G))}. \nonumber
\end{eqnarray}
Combining \eqref{5.11}--\eqref{5.14}, we arrive at
\begin{equation}\label{5.15}
\begin{aligned}
\lambda\dbE\int_{Q}e^{2\lambda\phi}z^{2}dxdt
&\leq C\Big(\|e^{\lambda\phi}(a_{4}z+Z)\|_{L_{\dbF}^{2}(0,T;L^{2}(G))}^{2}+
\|e^{\lambda\phi}(a_{5}\hat{z}+\h{Z})\|^{2}_{L_{\dbF}^{2}(0,T;H^{-1}(G))}\\
&\qq\q+\lambda(1+\lambda)e^{2\lambda e^{\mu_0(R_1^2-c_1 T^{2}/4)/2}}
\|z\|^{2}_{L^{2}_{\dbF}(J;L^{2}(G))}\\
&\qquad\quad+\|e^{\lambda\phi} z\|^{2}_{L^{2}_{\dbF}(0,T;L^{2}(G))}
+\lambda^{2}\dbE\int_{0}^{T}\int_{G_{0}}e^{2\lambda\phi}z^{2}dxdt\Big).
\end{aligned}
\end{equation}
From \eqref{5.1}, since $R_{1}^{2}-c_{1}T^{2}/4<0$, one can find $\lambda_{3}>0$ such that for any $\lambda>\lambda_{3}$, \eqref{5.7} holds.

\ss

{\bf Step 3.} By \eqref{5.3} we have
\begin{equation}\label{5.16}
\dbE\int_{Q}e^{2\lambda\phi}z^{2}dxdt\geq e^{2\lambda e^{\mu_{0}R_0^2/2}}\dbE\int_{T_{1}}^{T_{1}'}\int_{G}z^{2}dxdt.
\end{equation}
It follows from \eqref{4.3} that
\begin{equation}\label{5.17}
\begin{aligned}
\|z\|^{2}_{L^{2}_{\dbF}(J;L^{2}(G))}
\leq\dbE\int_{0}^{T}\cE(t)dt \leq C\big(\cE(T)+\|a_{4}z+Z\|^{2}_{L^{2}_{\dbF}(0,T;L^{2}(G))}
+\|a_{5}\hat{z}+\h{Z}\|^{2}_{L^{2}_{\dbF}(0,T;H^{-1}(G))}\big).
\end{aligned}
\end{equation}
Let $0<\vep_{3}<\vep_{1}$, and set
$$
T_{3}=\frac{T}{2}-\vep_{3}T,\q
T_{3}'=\frac{T}{2}+\vep_{3}T.
$$
By \eqref{4.2}, \eqref{4.3-1} and \eqref{5.16}, we obtain that
\begin{equation}\label{5.18}
\begin{aligned}
\cE(T)
&\leq C\Big(\int_{T_{3}'}^{T_{3}}\cE(t)dt+\|a_{4}z+Z\|^{2}_{L^{2}_{\dbF}(0,T;L^{2}(G))}
+\|a_{5}\hat{z}+\h{Z}\|^{2}_{L^{2}_{\dbF}(0,T;H^{-1}(G))}\Big)\\
&\leq C\Big(\dbE\int_{T_{0}}^{T_{0}'}\|z\|_{L^{2}(G)}^{2}dt+\|a_{4}z+Z\|^{2}_{L^{2}_{\dbF}(0,T;L^{2}(G))}
+\|a_{5}\hat{z}+\h{Z}\|^{2}_{L^{2}_{\dbF}(0,T;H^{-1}(G))}\Big)\\
&\leq C\Big(e^{-2\lambda e^{\mu_{0}R_0^2/2}}\dbE\int_{Q}e^{2\lambda\phi}z^{2}dxdt+\|a_{4}z+Z\|^{2}_{L^{2}_{\dbF}(0,T;L^{2}(G))}
+\|a_{5}\hat{z}+\h{Z}\|^{2}_{L^{2}_{\dbF}(0,T;H^{-1}(G))}\Big).
\end{aligned}
\end{equation}
Finally, combining \eqref{5.7}, \eqref{5.17}, and \eqref{5.18}, we have
\begin{eqnarray}\label{5.19} 
\cE(T)
&\3n\leq\3n&  C\Big[e^{-2\lambda e^{\mu_{0}R_0^2/2}}\Big(
\frac{1}{\lambda}\dbE\|z\|^{2}_{L^{2}(J\times G)}+\lambda\dbE\int_{0}^{T}\int_{G_{0}}e^{2\lambda\phi}z^{2}dxdt\Big)\nonumber\\
&&\qq\q+\|a_{4}z+Z\|^{2}_{L^{2}_{\dbF}(0,T;L^{2}(G))}
+\|a_{5}\hat{z}+\h{Z}\|^{2}_{L^{2}_{\dbF}(0,T;H^{-1}(G))}\Big]\\
&\3n\leq\3n&  C\Big[e^{-2\lambda e^{\mu_{0}R_0^2/2}}\Big(
\frac{1}{\lambda}\cE(T)+\lambda\dbE\int_{0}^{T}\int_{G_{0}}e^{2\lambda\phi}z^{2}dxdt\Big)\nonumber\\
&&\qq\q+\|a_{4}z+Z\|^{2}_{L^{2}_{\dbF}(0,T;L^{2}(G))}
+\|a_{5}\hat{z}+\h{Z}\|^{2}_{L^{2}_{\dbF}(0,T;H^{-1}(G))}\Big].\nonumber
\end{eqnarray}
Consequently, there exists a constant $\lambda_{4}>0$ such that for any $\lambda>\lambda_{4}$, \eqref{1.15} holds. This completes the proof. \endpf

\end{document}